\documentclass[11pt]{article}

\usepackage{latexsym}
\usepackage{amssymb}

\newtheorem{thm}{Theorem}[section]
\newtheorem{prop}{Proposition}[section]
\newtheorem{lem}{Lemma}[section]
\newtheorem{cor}{Corollary}[section]
\newtheorem{rmr}{Remark}[section]

\begin{document}
{
\begin{center}
{\Large\bf
On the generalized resolvents of isometric operators with gaps.}
\end{center}
\begin{center}
{\bf S.M. Zagorodnyuk}
\end{center}

\section{Introduction.}
We shall investigate generalized resolvents of an isometric operator.
Let $V$ be a closed isometric operator in a (separable) Hilbert space $H$. There always exists
(at least one) unitary operator $U\supseteq V$ in a Hilbert space $\widetilde H\supseteq H$.
Recall that the following operator-valued function $\mathbf{R}_\zeta$:
$$ \mathbf{R}_\zeta h = P^{\widetilde{H}}_H \left( E_{\widetilde{H}} - \zeta U \right)^{-1} h,\quad
h\in H, $$
is said to be the {\bf generalized resolvent} of the isometric operator $V$ (corresponding to the
extension $U$). Let $\{ F_t \}_{t\in [0,2\pi]}$ be the left-continuous orthogonal resolution of
unity of $U$. Then the following operator-valued function:
$$ \mathbf{F}_t = P^{\widetilde{H}}_H F_t,\qquad t\in [0,2\pi], $$
is said to be a (left-continuous) {\bf spectral function } of the isometric operator $V$ (corresponding to the
extension $U$).
Let $F(\delta)$, $\delta\in\mathfrak{B}(\mathbb{T})$, be the orthogonal spectral measure of $U$.
Then
$$ \mathbf{F}(\delta) = P^{\widetilde{H}}_H F(\delta),\qquad \delta\in \mathfrak{B}(\mathbb{T}), $$
is said to be a {\bf spectral measure } of the isometric operator $V$ (corresponding to the
extension $U$). Of course, we have
$$ \mathbf{F}(\delta_t) = \mathbf{F}_t,\qquad \delta_t=\{ z=e^{i\varphi}:\ 0\leq \varphi < t \},\quad
t\in [0,2\pi], $$
what follows from the analogous property of the orthogonal measures.
We notice that there exists a one-to-one correspondence between spectral functions (spectral measures)
and generalized resolvents:
\begin{equation}
\label{f1_m1}
(\mathbf R_z h,g)_H = \int_{\mathbb{T}} \frac{1}{1-z\zeta} d(\mathbf{F}(\cdot) h,g)_H =
\int_0^{2\pi} \frac{1}{1-ze^{it}} d(\mathbf{F}_t h,g)_H,\quad \forall h,g\in H,
\end{equation}
according to the inversion formula~\cite[p.50]{cit_3000_AK}.

Let $H_1$ and $H_2$ be two arbitrary subspaces of the Hilbert space $H$.
By $\mathcal{S}(H_1;H_2)$ we denote the set of all analytic in $\mathbb{D}=\{ z\in \mathbb{C}:\ |z|<1 \}$
operator-valued functions $F(\zeta)$ which values are linear contractions with the domain $D(F(\zeta)) = H_1$ and
with the range $R(F(\zeta))\subseteq H_2$, $\forall\zeta\in \mathbb{D}$.

\noindent
Chumakin's formula~\cite[Theorem 3]{cit_4000_Ch}:
\begin{equation}
\label{f1_1}
\mathbf R_{\zeta} = \left[
E_H - \zeta ( V \oplus F(\zeta) )
\right]^{-1},\qquad
\zeta\in \mathbb{D},
\end{equation}
establishes a one-to-one correspondence between all generalized resolvents of $V$ and all functions
$F(\zeta)$ from the set $\mathcal{S}(H\ominus D(V);H\ominus R(V))$.

Set
$$ M_\zeta = M_\zeta(V) = (E_H - \zeta V) D(V),\quad N_\zeta = N_\zeta(V) =
H\ominus M_\zeta,\quad \zeta\in \mathbb{C}; $$
$$ M_\infty = M_\infty(V) = R(V),\quad N_\infty = N_\infty(V) = H\ominus R(V). $$
Consider the following operator
\begin{equation}
\label{f1_2}
V_z = (V-\overline{z}E_H) (E_H - zV)^{-1},\qquad z\in \mathbb{D}.
\end{equation}
Notice that $D(V_z) = M_z$ and $R(V_z) = M_{\frac{1}{ \overline{z} }}$.
It is straightforward to check that $V_z$ is isometric and
\begin{equation}
\label{f1_3}
V = (V_z + \overline{z}E_H) (E_H + zV_z)^{-1} = \left( V_z \right)_{-z}.
\end{equation}
Moreover, if $V$ is unitary, then $V_z$ is unitary, and vice versa (by~(\ref{f1_3})).

Let $\widehat V_z\supseteq V_z$ be a unitary operator in a Hilbert space $\widehat H\supseteq H$.
Then we may define the operator
\begin{equation}
\label{f1_4}
\widehat V = (\widehat V_z + \overline{z} E_{\widehat H}) (E_{\widehat H} + z\widehat V_z)^{-1},
\end{equation}
which is a unitary extension of $V$.
Formula~(\ref{f1_4}) establishes a one-to-one correspondence between all unitary extensions $\widehat V_z$ of
$V_z$ in
a Hilbert space $\widehat H\supseteq H$, and
all unitary extensions $\widehat V$ of $V$ in a Hilbert space $\widehat H$.

Let us fix an arbitrary point $z_0\in \mathbb{D}$. Let $C$ be an arbitrary linear bounded operator with the
domain $D(C) = N_{z_0}$ and the range $R(C)\subseteq N_{\frac{1}{ \overline{z_0} }}$.
Set
\begin{equation}
\label{f1_5}
V^+_{z_0;C} = V_{z_0} \oplus C;
\end{equation}
\begin{equation}
\label{f1_6}
V_{C}=V_{C;z_0} =  (V^+_{z_0;C} + \overline{z_0} E_{H}) (E_{H} + z_0 V^+_{z_0;C})^{-1}.
\end{equation}
If $z_0\not=0$, we may write:
\begin{equation}
\label{f1_6_1}
V_{C}=V_{C;z_0} =  \frac{1}{z_0} E_H + \frac{|z_0|^2-1}{z_0} (E_{H} + z_0 V^+_{z_0;C})^{-1};
\end{equation}
\begin{equation}
\label{f1_6_2}
V^+_{z_0;C} =  -\frac{1}{z_0} E_H + \frac{ 1 - |z_0|^2}{z_0} (E_{H} - z_0 V_{C;z_0})^{-1}.
\end{equation}
Recall that the operator $V_C$ is said to be an {\bf orthogonal extension} of $V$ defined by
the operator $C$.

\noindent
Inin's formula~\cite[Theorem]{cit_5000_I}:
\begin{equation}
\label{f1_7}
\mathbf R_{\zeta} = \left[
E - \zeta V_{C(\zeta;z_0)}
\right]^{-1},\qquad
\zeta\in \mathbb{D},
\end{equation}
establishes a one-to-one correspondence between all generalized resolvents of $V$ and all functions
$C(\zeta)=C(\zeta;z_0)$ from the set $\mathcal{S}(N_{z_0};N_{\frac{1}{ \overline{z_0} }})$.
Observe that in the case $z_0=0$ it coincides with Chumakin's formula.

We shall show that Inin's formula can be derived directly from Chumakin's formula.
Then we shall obtain an analog of some McKelvey's
results~\cite[Theorem 2.1 (A),(B); Theorem 3.1 (A),(B); Remark 2.2]{cit_5500_M}, see also~\cite{cit_5700_VL}.
Also we obtain an auxiliary proposition which
uses some constructions of L.A.~Shtraus in~\cite[Lemma]{cit_5900_S_LA}.
All that will be used to obtain some slight correction and generalization
of Ryabtseva's results about generalized resolvents of an isometric operator with a gap
in~\cite{cit_5950_R}. Here we used some
ideas of Varlamova-Luks for the case of Hermitian operators with a
gap~\cite{cit_5700_VL},\cite{cit_5960_VL},\cite{cit_5970_VL}.

{\bf Notations. }
As usual, we denote by $\mathbb{R}, \mathbb{C}, \mathbb{N}, \mathbb{Z}, \mathbb{Z}_+$,
the sets of real numbers, complex numbers, positive integers, integers and non-negative integers,
respectively; $\mathbb{D} = \{ z\in \mathbb{C}:\ |z|<1 \}$, $\mathbb{T} = \{ z\in \mathbb{C}:\ |z|=1 \}$,
$\mathbb{D}_e = \{ z\in \mathbb{C}:\ |z|>1 \}$, $\mathbb{T}_e = \{ z\in \mathbb{C}:\ |z|\not=1 \}$.
By $\mathfrak{B}(\mathbb{T})$ we denote the set of all Borel subsets of $\mathbb{T}$.

All Hilbert spaces in this paper are assumed to be separable.
If H is a Hilbert space then $(\cdot,\cdot)_H$ and $\| \cdot \|_H$ mean
the scalar product and the norm in $H$, respectively.
Indices may be omitted in obvious cases.
For a linear operator $A$ in $H$, we denote by $D(A)$
its  domain, by $R(A)$ its range, by $\mathop{\rm Ker}\nolimits A$
its null subspace (kernel), and $A^*$ means the adjoint operator
if it exists. If $A$ is invertible then $A^{-1}$ means its
inverse. $\overline{A}$ means the closure of the operator, if the
operator is closable. If $A$ is bounded then $\| A \|$ denotes its
norm. The set of all points of the regular type of $A$ is denoted by $\mathcal{M}_r(A)$.
For a set $M\subseteq H$
we denote by $\overline{M}$ the closure of $M$ in the norm of $H$.
By $A|_M$ we denote the restriction of the operator $A$ to $M$.
For an arbitrary set of elements $\{ x_n \}_{n\in I}$ in
$H$, we denote by $\mathop{\rm Lin}\nolimits\{ x_n \}_{n\in I}$
the set of all linear combinations of elements $x_n$,
and $\mathop{\rm span}\nolimits\{ x_n \}_{n\in I}
:= \overline{ \mathop{\rm Lin}\nolimits\{ x_n \}_{n\in I} }$.
Here $I$ is an arbitrary set of indices.
By $E_H$ we denote the identity operator in $H$, i.e. $E_H x = x$,
$x\in H$. In obvious cases we may omit the index $H$. If $H_1$ is a subspace of $H$, then $P_{H_1} =
P_{H_1}^{H}$ is an operator of the orthogonal projection on $H_1$
in $H$. By $w.-\lim$ and $u.-\lim$ we denote the limits in the weak and the uniform operator
topologies, respectively.

\section{A connection between Chumakin's formula and Inin's formula.}
The following proposition holds, see~\cite[p.34]{cit_5000_I}.
\begin{prop}
\label{p2_1}
Let $V$ be a closed isometric operator in a Hilbert space $H$. Let $z_0\in \mathbb{D}$ be fixed.
For an arbitrary point $\zeta\in \mathbb{C}\backslash\{ 0 \}$, $\zeta\not= z_0$,
the following two conditions are equivalent:
\begin{itemize}
\item[(i)]  $\zeta^{-1}\in \mathcal{M}_r(V)$;
\item[(ii)] $\frac{1-\zeta\overline{z_0}}{\zeta-z_0}\in \mathcal{M}_r(V_{z_0})$.
\end{itemize}
\end{prop}
{\bf Proof.}
$(i)\Rightarrow(ii)$.
We may write
$$ V_{z_0} - \frac{1-\zeta\overline{z_0}}{\zeta-z_0} E_H =
(V- \overline{z_0} E_H)(E_H - z_0 V)^{-1} - \frac{1-\zeta\overline{z_0}}{\zeta-z_0} (E_H - z_0 V)(E_H - z_0 V)^{-1} $$
$$ = \frac{\zeta(1-|z_0|^2)}{\zeta - z_0} (V-\frac{1}{\zeta} E_H)(E_H - z_0 V)^{-1}. $$
The operator on the right-hand side has a bounded inverse defined on $(V-\zeta^{-1} E_H)D(V)$.

\noindent
$(ii)\Rightarrow(i)$. We write:
$$ V - \frac{1}{\zeta} E_H =
(V_{z_0} + \overline{z_0}E_H) (E_H + z_0 V_{z_0})^{-1} - \frac{1}{\zeta} (E_H + z_0 V_{z_0})
(E_H + z_0 V_{z_0})^{-1} $$
$$ = \frac{\zeta - z_0}{\zeta} \left(
V_{z_0} - \frac{1-\zeta\overline{z_0}}{\zeta-z_0} E_H
\right)
(E_H + z_0 V_{z_0})^{-1}, $$
and the operator on the right-hand side has a bounded inverse which is defined
on $(V_{z_0} - \frac{1-\zeta\overline{z_0}}{\zeta-z_0} E_H)D(V_{z_0})$.
$\Box$

Let $V$ be a closed isometric operator in a Hilbert space $H$, and
$z_0\in \mathbb{D}\backslash\{ 0 \}$ be a fixed point.
Consider the following linear fractional transformation:
\begin{equation}
\label{f2_1}
t = t(u) = \frac{u-\overline{z_0}}{1-z_0 u},
\end{equation}
which maps $\mathbb{T}$ on $\mathbb{T}$, and $\mathbb{D}$ on $\mathbb{D}$.

\noindent
Let $\widehat V_{z_0}$ be an arbitrary unitary extension of $V_{z_0}$ in a Hilbert space
$\widehat H\supseteq H$, and $\widehat V$ be the corresponding unitary extension of $V$ defined
by relation~(\ref{f1_4}).
Choose an arbitrary $u\in \mathbb{T}_e\backslash\{ 0, \overline{z_0}, \frac{1}{z_0} \}$.
Then
$t = t(u)\in \mathbb{T}_e\backslash\{ 0, -\overline{z_0}, -\frac{1}{z_0} \}$.
Moreover
\begin{equation}
\label{f2_1_1}
u\in \mathbb{T}_e\backslash\{ 0, \overline{z_0}, \frac{1}{z_0} \}
\Leftrightarrow
t \in \mathbb{T}_e\backslash\{ 0, -\overline{z_0}, -\frac{1}{z_0} \}.
\end{equation}
We may write:
$$ ( \widehat V_{z_0} - tE_{\widehat H} )^{-1} =
\left(
(\widehat V-\overline{z_0}E_{\widehat H} ) (E_{\widehat H} - z_0 \widehat V)^{-1}
 \right. $$
$$ \left.
- \frac{u-\overline{z_0}}{1-z_0 u} (E_{\widehat H} - z_0\widehat V) (E_{\widehat H} - z_0\widehat V)^{-1}
\right)^{-1}
$$
$$ =
\left(
\frac{1-|z_0|^2}{1-z_0 u}
(\widehat V - uE_{\widehat H})
(E_{\widehat H} - z_0\widehat V)^{-1}
\right)^{-1}.
$$
$$ = \frac{1-z_0 u}{1-|z_0|^2} (E_{\widehat H} - z_0\widehat V)
(\widehat V - uE_{\widehat H})^{-1} $$
$$ = - \frac{z_0(1-z_0 u)}{1-|z_0|^2} E_{\widehat H} +
\frac{(1-z_0 u)^2}{1-|z_0|^2} (\widehat V - uE_{\widehat H})^{-1}. $$
Therefore
$$ -\frac{1}{t} (E_{\widehat H} - \frac{1}{t}\widehat V_{z_0})^{-1}
=
- \frac{z_0(1-z_0 u)}{1-|z_0|^2} E_{\widehat H} -
\frac{(1-z_0 u)^2}{u(1-|z_0|^2)} (E_{\widehat H} - \frac{1}{u}\widehat V)^{-1}; $$
$$ (E_{\widehat H} - \frac{1}{u}\widehat V)^{-1} =
-\frac{z_0 u}{1-z_0 u} E_{\widehat H} + \frac{u(1-|z_0|^2)}{(1-z_0 u)^2 t}
(E_{\widehat H} - \frac{1}{t}\widehat V_{z_0})^{-1} $$
$$ = -\frac{z_0 u}{1-z_0 u} E_{\widehat H} + \frac{u(1-|z_0|^2)}{(1-z_0 u)(u-\overline{z_0})}
(E_{\widehat H} - \frac{1}{t}\widehat V_{z_0})^{-1}. $$
Set $\widetilde u = \frac{1}{u}$, $\widetilde t = \frac{1}{t}$. Observe that
$\widetilde u\in \mathbb{T}_e\backslash\{ 0, \frac{1}{\overline{z_0}}, z_0 \}$,
$\widetilde t\in \mathbb{T}_e\backslash\{ 0, -z_0, -\frac{1}{\overline{z_0}} \}$.
Moreover
\begin{equation}
\label{f2_1_2}
\widetilde u\in \mathbb{T}_e\backslash\{ 0, \frac{1}{\overline{z_0}}, z_0 \}
\Leftrightarrow
\widetilde t\in \mathbb{T}_e\backslash\{ 0, -z_0, -\frac{1}{\overline{z_0}} \},
\end{equation}
and these conditions are equivalent to conditions from relation~(\ref{f2_1_1}).
Then
$$ (E_{\widehat H} - \widetilde{u}\widehat V)^{-1} =
-\frac{z_0}{\widetilde{u}-z_0} E_{\widehat H} + \frac{\widetilde{u} (1-|z_0|^2)}{ (\widetilde{u} - z_0)
(1 - \overline{z_0} \widetilde{u}) }
(E_{\widehat H} - \widetilde{t}\widehat V_{z_0})^{-1}. $$
By applying the projection operator $P^{\widehat H}_H$ to the both sides of the last relation,
we obtain the following relation:
\begin{equation}
\label{f2_2}
\mathbf{R}_{\widetilde u}(V) = -\frac{z_0}{\widetilde{u}-z_0} E_{H}
+ \frac{\widetilde{u} (1-|z_0|^2)}{ (\widetilde{u} - z_0)(1 - \overline{z_0} \widetilde{u}) }
\mathbf{R}_{\frac{ \widetilde{u} - z_0 }{ 1 - \overline{z_0} \widetilde{u} }}(V_{z_0}),\quad
\widetilde u\in \mathbb{T}_e\backslash\{ 0, \frac{1}{\overline{z_0}}, z_0 \},
\end{equation}
where $\mathbf{R}_{\widetilde u}(V)$, $\mathbf{R}_{\widetilde t}(V_{z_0})$, are the generalized resolvents of
the operators $V$, $V_{z_0}$, respectively.

\noindent
Since $\mathbf{R}_{\widetilde u}(V)$ is analytic in $\mathbb{T}_e$, it is uniquely defined by
the generalized resolvent
$\mathbf{R}_{\widetilde{t}}(V_{z_0})$, by relation~(\ref{f2_2}).
By the same relation~(\ref{f2_2}), the generalized resolvent $\mathbf{R}_{\widetilde{t}}(V_{z_0})$
is uniquely defined by the generalized resolvent $\mathbf{R}_{\widetilde u}(V)$.

\noindent
Thus, relation~(\ref{f2_2}) establishes a one-to-one correspondence between all generalized
resolvents of $V_{z_0}$, and all generalized resolvents of $V$.

Let us apply Chumakin's formula~(\ref{f1_1}) to the operator $V_{z_0}$:
\begin{equation}
\label{f2_3}
\mathbf{R}_{\widetilde t}(V_{z_0}) =
\left[
E_H - \widetilde t ( V_{z_0} \oplus F(\widetilde t) )
\right]^{-1},\qquad
\widetilde t\in \mathbb{D},
\end{equation}
where $F(\widetilde t)$ belongs to the set $\mathcal{S}(N_{z_0};N_{\frac{1}{ \overline{z_0} }})$.
Let us restrict relation~(\ref{f2_2}) to $\widetilde u\in \mathbb{D}\backslash\{ 0,z_0 \}$ what is
equivalent to the condition $\widetilde t\in \mathbb{D}\backslash\{ 0, -z_0 \}$. In this case it
also establishes the above-mentioned one-to-one correspondence. By~(\ref{f2_2}),(\ref{f2_3})
we get
$$ \mathbf{R}_{\widetilde u}(V) = -\frac{z_0}{\widetilde{u}-z_0} E_{H} $$
\begin{equation}
\label{f2_4}
+ \frac{\widetilde{u} (1-|z_0|^2)}{ (\widetilde{u} - z_0)(1 - \overline{z_0} \widetilde{u}) }
\left[
E_H - \frac{ \widetilde{u} - z_0 }{ 1 - \overline{z_0} \widetilde{u} }
\left( V_{z_0} \oplus
F\left(
\frac{ \widetilde{u} - z_0 }{ 1 - \overline{z_0} \widetilde{u} }
\right)
\right)
\right]^{-1},\
\widetilde u\in \mathbb{D}\backslash\{ 0,z_0 \},
\end{equation}
where $F(\widetilde t)\in\mathcal{S}(N_{z_0};N_{\frac{1}{ \overline{z_0} }})$.
Relation~(\ref{f2_4}) establishes a one-to-one correspondence between all functions
$F(\widetilde t)$ from the set $\mathcal{S}(N_{z_0};N_{\frac{1}{ \overline{z_0} }})$, and
all generalized resolvents of the operator $V$.

\noindent
Set $C(\widetilde u) = F(
\frac{ \widetilde{u} - z_0 }{ 1 - \overline{z_0} \widetilde{u} })$, $u\in \mathbb{D}$.
Observe that $C(\widetilde u)\in\mathcal{S}(N_{z_0};N_{\frac{1}{ \overline{z_0} }})$.
We may write
$$ E_H - \frac{ \widetilde{u} - z_0 }{ 1 - \overline{z_0} \widetilde{u} }
\left( V_{z_0} \oplus
F\left(
\frac{ \widetilde{u} - z_0 }{ 1 - \overline{z_0} \widetilde{u} }
\right)
\right) =
E_H - \frac{ \widetilde{u} - z_0 }{ 1 - \overline{z_0} \widetilde{u} }
V^+_{z_0;C(\widetilde u)} $$
$$ = (E_H - z_0 V_{C(\widetilde u);z_0}) (E_H - z_0 V_{C(\widetilde u);z_0})^{-1} $$
$$ - \frac{ \widetilde{u} - z_0 }{ 1 - \overline{z_0} \widetilde{u} }
(V_{C(\widetilde u);z_0} - \overline{z_0} E_H) (E_H - z_0 V_{C(\widetilde u);z_0})^{-1} $$
$$ = \frac{1-|z_0|^2}{ 1 - \overline{z_0}\widetilde u} (E_H - \widetilde u V_{C(\widetilde u);z_0})
(E_H - z_0 V_{C(\widetilde u);z_0})^{-1} $$
By substitution the last relation into relation~(\ref{f2_4}) and after elementary calculations we get:
\begin{equation}
\label{f2_5}
\mathbf{R}_{\widetilde u}(V) =
(E_H - \widetilde u V_{C(\widetilde u);z_0})^{-1},\qquad
\widetilde u\in \mathbb{D}\backslash\{ 0,z_0 \}.
\end{equation}
Of course, for the case $\widetilde u = 0$ relation~(\ref{f2_5}) is also true.
It remains to check relation~(\ref{f2_5}) for the case $\widetilde u = z_0$, to obtain
Inin's formula.

\noindent
By Chumakin's formula for $\mathbf{R}_{\widetilde u}(V)$ we see that
$(\mathbf{R}_{\widetilde u}(V))^{-1}$ is an analytic operator-valued function in $\mathbb{D}$.
By~(\ref{f2_5}) we may write
\begin{equation}
\label{f2_5_1}
(\mathbf{R}_{z_0}(V))^{-1} = u.-\lim_{\widetilde u\to z_0} (\mathbf{R}_{\widetilde u}(V))^{-1}
= E_H - u.-\lim_{\widetilde u\to z_0} \widetilde u V_{C(\widetilde u);z_0},
\end{equation}
where the limits are understood in the uniform operator topology.

\noindent
The operator-valued function
$V^+_{z_0;C(\widetilde u)} = V_{z_0}\oplus C(\widetilde u)$ is analytic in $\mathbb{D}$,
and its values are contractions in $H$. Then
$$ \| (E_H + z_0 V^+_{z_0;C(\widetilde u)}) h \| \geq
| \| h \| - |z_0| \| V^+_{z_0;C(\widetilde u)}) h \| | \geq
(1 - |z_0|) \| h \|,\qquad h\in H; $$
\begin{equation}
\label{f2_6}
\| (E_H + z_0 V^+_{z_0;C(\widetilde u)})^{-1} \| \leq \frac{1}{1 - |z_0|},\qquad
\widetilde u\in \mathbb{D}.
\end{equation}
Using~\cite[Footnote on page 83]{cit_7000_S} we obtain that the function
$(E_H + z_0 V^+_{z_0;C(\widetilde u)})^{-1}$ is analytic in $\mathbb{D}$.
Therefore
$V_{C(\widetilde u);z_0} = (V^+_{z_0;C(\widetilde u)} + \overline{z_0} E_H)
(E_H + z_0 V^+_{z_0;C(\widetilde u)})^{-1}$
is analytic in $\mathbb{D}$, as well.
Passing to the limit in relation~(\ref{f2_5_1}) we get
$$ (\mathbf{R}_{z_0}(V))^{-1} =
E_H - z_0 V_{C(z_0);z_0}. $$
Therefore relation~(\ref{f2_5}) holds for the case $\widetilde u = z_0$, and we proved
Inin's formula.

\section{An analog of McKelvey's results.}
The following proposition is an analog of Theorem~3.1 (A),(B) in~\cite{cit_5500_M},
see also~\cite[Lemma 1.1]{cit_5700_VL}.
\begin{prop}
\label{p3_1}
Let $V$ be a closed isometric operator in a Hilbert space $H$, and $\mathbf{F}(\delta)$,
$\delta\in \mathfrak{B}(\mathbb{T})$,
be its spectral measure.
The following two conditions are equivalent:
\begin{itemize}
\item[(i)]  $\mathbf{F}(\Delta) = 0$, for some open arc $\Delta$ of $\mathbb{T}$;
\item[(ii)] The generalized resolvent $\mathbf{R}_z(V)$, corresponding to the spectral measure
$\mathbf{F}(\delta)$, has an analytic continuation to the set $\mathbb{D}\cup\mathbb{D}_e\cup \overline{\Delta}$,
where $\overline{\Delta} = \{ z\in \mathbb{C}:\ \overline{z}\in \Delta \}$, for some open arc
$\Delta$ of $\mathbb{T}$.
\end{itemize}
\end{prop}
{\bf Proof. }
(i)$\Rightarrow$(ii).
In this case relation~(\ref{f1_m1}) takes the following form:
\begin{equation}
\label{f3_1}
(\mathbf R_z h,g)_H = \int_{\mathbb{T}\backslash\Delta} \frac{1}{1-z\zeta} d(\mathbf{F}(\cdot) h,g)_H,\quad
\forall h,g\in H.
\end{equation}
Choose an arbitrary $z_0\in \overline{\Delta}$. Since $\frac{1}{1-z_0\zeta}$ is bounded and continuous
on $\mathbb{T}\backslash\Delta$, there exists an integral
$$ I_{z_0}(h,g) := \int_{\mathbb{T}\backslash\Delta} \frac{1}{1-z_0\zeta} d(\mathbf{F}(\cdot) h,g)_H. $$
Then
$$ | (\mathbf R_z h,h)_H - I_{z_0}(h,h) | =
|z-z_0| \left|
\int_{\mathbb{T}\backslash\Delta} \frac{\zeta}{(1-z\zeta)(1-z_0\zeta)} d(\mathbf{F}(\cdot) h,h)_H
\right|
$$
$$ \leq
|z-z_0| \int_{\mathbb{T}\backslash\Delta} \frac{|\zeta|}{|1-z\zeta||1-z_0\zeta|} d(\mathbf{F}(\cdot) h,h)_H,\quad
z\in \mathbb{T}_e. $$
There exists a neighborhood $U(z_0)$ of $z_0$ such that $|z-\overline{\zeta}|\geq M_1 > 0$, $\forall\zeta
\in \mathbb{T}\backslash\Delta$, $\forall z\in U(z_0)$. Therefore the integral in the last relation
is bounded in $U(z_0)$.
Thus, we obtain that
$$ (\mathbf R_z h,h)_H \rightarrow I_{z_0}(h,h),\quad z\in\mathbb{T}_e,\ z\rightarrow z_0,\quad \forall h\in H. $$
Using properties of sesquilinear forms we get
$$ (\mathbf R_z h,g)_H \rightarrow I_{z_0}(h,g),\quad z\in\mathbb{T}_e,\ z\rightarrow z_0,\quad \forall h,g\in H. $$
Set
$$ \mathbf R_{\widetilde z} :=
w.-\lim_{z\in \mathbb{T}_e,\ z\to \widetilde z} \mathbf R_z,\quad \forall \widetilde z\in\overline{\Delta}, $$
where the limit is understood in the weak operator topology.
We may write
$$ \left(
\frac{1}{z-z_0}(\mathbf R_z -  \mathbf R_{z_0})h,h
\right)_H =
\int_{\mathbb{T}\backslash\Delta} \frac{\zeta}{(1-z\zeta)(1-z_0\zeta)} d(\mathbf{F}(\cdot) h,h)_H, $$
$$ z\in U(z_0),\ h\in H. $$
The function under the integral is bounded in $U(z_0)$, and it tends to
$\frac{\zeta}{(1-z_0\zeta)^2}$. By the Lebesgue convergence theorem we deduce that
$$ \lim_{z\to z_0} \left(
\frac{1}{z-z_0}(\mathbf R_z -  \mathbf R_{z_0})h,h
\right)_H =
\int_{\mathbb{T}\backslash\Delta} \frac{\zeta}{(1-z_0\zeta)^2} d(\mathbf{F}(\cdot) h,h)_H; $$
and therefore
$$ \lim_{z\to z_0} \left(
\frac{1}{z-z_0}(\mathbf R_z -  \mathbf R_{z_0})h,g
\right)_H =
\int_{\mathbb{T}\backslash\Delta} \frac{\zeta}{(1-z_0\zeta)^2} d(\mathbf{F}(\cdot) h,g)_H, $$
for $h,g\in H$.
Consequently, there exists the derivative of $\mathbf R_z$ at $z=z_0$.

\noindent
(ii)$\Rightarrow$(i).
Choose an arbitrary $h\in H$, and consider the function $\sigma_h(t) := (\mathbf{F}_t h,h)_H$,
$t\in [0,2\pi)$, where $\mathbf{F}_t$ is the left-continuous spectral function of $V$,
corresponding to the spectral measure $\mathbf{F}(\delta)$. Also consider the following function:
$$ f_h(z) = \frac{1}{2} \int_0^{2\pi} \frac{1+e^{it} z}{1-e^{it} z} d\sigma_h(t)
= \int_0^{2\pi} \frac{1}{1-e^{it} z} d\sigma_h(t) - \frac{1}{2}
\int_0^{2\pi} d\sigma_h(t) $$
\begin{equation}
\label{f3_1_1}
= \int_0^{2\pi} \frac{1}{1-e^{it} z} d\sigma_h(t) - \frac{1}{2} \| h \|^2_H
= (\mathbf{R}_z h,h)_H - \frac{1}{2} \| h \|^2_H.
\end{equation}
Choose arbitrary numbers $t_1,t_2$, $0\leq t_1<t_2\leq 2\pi$, such that
\begin{equation}
\label{f3_2}
l=l(t_1,t_2)= \{ z=e^{it}:\ t_1\leq t\leq t_2 \} \subset \Delta.
\end{equation}
Suppose additionally that $t_1$ and $t_2$ are points of continuity of the function $\mathbf{F}_t$.
By the inversion formula~\cite[p.50]{cit_3000_AK} we may write:
$$ \sigma_h(t_2) - \sigma_h(t_1) = \lim_{r\to 1-0} \int_{t_1}^{t_2} {\mathrm Re}
\left\{ f_h(re^{-i\tau}) \right\} d\tau. $$
Observe that
$$ {\mathrm Re} \left\{ f_h(re^{-i\tau}) \right\} =
{\mathrm Re} \left\{ (\mathbf{R}_{re^{-i\tau}} h,h)_H \right\} - \frac{1}{2} \| h \|^2_H $$
\begin{equation}
\label{f3_2_1}
= \frac{1}{2} \left(
((\mathbf{R}_{re^{-i\tau}} + \mathbf{R}_{re^{-i\tau}}^*) h,h)_H
\right) - \frac{1}{2} \| h \|^2_H,\quad t_1\leq \tau \leq t_2.
\end{equation}
By~(\ref{f3_2}) we see that $e^{-i\tau}$ belongs to $\overline\Delta$, for $t_1\leq \tau \leq t_2$.
Therefore
\begin{equation}
\label{f3_2_2}
\lim_{r\to 1-0} ((\mathbf{R}_{re^{-i\tau}} + \mathbf{R}_{re^{-i\tau}}^*) h,h)_H =
((\mathbf{R}_{e^{-i\tau}} + \mathbf{R}_{e^{-i\tau}}^*) h,h)_H.
\end{equation}
The generalized resolvents have the following property~\cite{cit_4000_Ch}:
\begin{equation}
\label{f3_3}
\mathbf{R}_z^* = E_H - \mathbf{R}_{\frac{1}{ \overline{z} }},\quad z\in \mathbb{T}_e.
\end{equation}
Passing to the limit in~(\ref{f3_3}) as $z$ tends to $e^{-i\tau}$, we get
\begin{equation}
\label{f3_4}
\mathbf{R}_{e^{-i\tau}}^* = E_H - \mathbf{R}_{e^{-i\tau}},\quad t_1\leq \tau \leq t_2.
\end{equation}
By~(\ref{f3_2_1}),(\ref{f3_2_2}) and~(\ref{f3_4}) we obtain that
\begin{equation}
\label{f3_5}
\lim_{r\to 1-0} {\mathrm Re} \left\{ f_h(re^{-i\tau}) \right\} = 0,\quad t_1\leq \tau \leq t_2.
\end{equation}
Consider the following sector:
$$ L(t_1,t_2) = \{ z=re^{-it}:\ t_1\leq t\leq t_2,\ 0\leq r \leq 1 \}. $$
The generalized resolvent is analytic at any point of the closed sector $L(t_1,t_2)$. Therefore
${\mathrm Re} (\mathbf{R}_z h,h)$ is continuous and bounded in $L(t_1,t_2)$.
By the Lebesgue convergence theorem we conclude that $\sigma_h(t_1)=\sigma_h(t_2)$.
If $1\notin\Delta$ we easily get the required result.
In the case $1\in\Delta$, we write $\Delta = \Delta_1\cup\{ 1 \}\cup\Delta_2$ where open arcs
$\Delta_1$ and $\Delta_2$ do not contain $1$. Then $\sigma_h(t)$ is constant in the intervals
corresponding to $\Delta_1$ and $\Delta_2$.
Suppose that there exists a non-zero jump of $\sigma_h(t)$ at $t=0$.
By~(\ref{f1_m1}) we may write:
$$ (\mathbf{R}_z h,h)_H  =
\int_0^{2\pi} \frac{1}{1-e^{it} z} d\sigma_h(t) =
\int_0^{2\pi} \frac{1}{1-e^{it} z} d\widehat\sigma_h(t) + \frac{1}{1-z} a,\quad a>0, $$
where $\widehat\sigma_h(t) = \sigma_h(t) + \sigma_h(+0) - \sigma_h(0)$, $t\in [0,2\pi]$.
In a neighborhood  of $1$ the left-hand side and the first summand of the right-hand side are bounded.
We obtained a contradiction.
$\Box$

The following theorem is an analog of Theorem~2.1 (A) in~\cite{cit_5500_M}.
\begin{thm}
\label{t3_1}
Let $V$ be a closed isometric operator in a Hilbert space $H$, and $\mathbf{R}_z(V)$ be an arbitrary
generalized resolvent of $V$. Let $\{ \lambda_k \}_{k=1}^\infty$ be a sequence of points of $\mathbb{D}$, such that
$\lambda_k\rightarrow \widehat\lambda$, as $k\rightarrow\infty$; $\widehat\lambda\in \mathbb{T}$.
Suppose that for some $z_0\in \mathbb{D}\backslash\{ 0 \}$,
the function $C(\lambda;z_0)$, corresponding to $\mathbf{R}_z(V)$
by Inin's formula~(\ref{f1_7}), satisfies the following relation:
\begin{equation}
\label{f3_6}
\exists u.-\lim_{k\to\infty} C(\lambda_k;z_0) =: C(\widehat\lambda;z_0).
\end{equation}
Then for arbitrary $z_0'\in \mathbb{D}\backslash\{ 0 \}$, the function $C(\lambda;z_0')$, corresponding to
$\mathbf{R}_z(V)$ by Inin's formula~(\ref{f1_7}), is such that
\begin{equation}
\label{f3_7}
\exists u.-\lim_{k\to\infty} C(\lambda_k;z_0') =: C(\widehat\lambda;z_0').
\end{equation}
In this case $C(\widehat\lambda;z_0')$ is a linear contraction which maps $N_{z_0'}$ into
$N_{\frac{1}{ \overline{z_0'} }}$, and the corresponding orthogonal extension
$V_{C(\widehat\lambda;z_0');z_0'}$ does not depend on the choice of $z_0'\in \mathbb{D}\backslash\{ 0 \}$.
\end{thm}
{\bf Proof. }
Suppose that relation~(\ref{f3_6}) holds for some $z_0\in \mathbb{D}\backslash\{ 0 \}$.
Choose an arbitrary $z_0'\in \mathbb{D}\backslash\{ 0 \}$.
Comparing Inin's formula for the choices $z_0$ and $z_0'$ we conclude that
\begin{equation}
\label{f3_8}
V_{C(\lambda;z_0);z_0} = V_{C(\lambda;z_0');z_0'},\qquad \lambda\in \mathbb{D}.
\end{equation}
By~(\ref{f1_6_1}) we write
\begin{equation}
\label{f3_9}
V_{C(\lambda;z_0);z_0} = \frac{1}{z_0} E_H +
\frac{|z_0|^2-1}{z_0}
\left(
E_H + z_0 V^+_{z_0;C(\lambda;z_0)}
\right)^{-1},\qquad \lambda\in \mathbb{D}.
\end{equation}
By substitution in~(\ref{f3_8}) such expressions for $z_0$ and $z_0'$, and multiplying by $z_0 z_0'$
we get
$$ z_0' E_H + z_0'(|z_0|^2-1)
\left(
E_H + z_0 V^+_{z_0;C(\lambda;z_0)}
\right)^{-1} $$
$$ =
z_0 E_H + z_0(|z_0'|^2-1)
\left(
E_H + z_0' V^+_{z_0';C(\lambda;z_0')}
\right)^{-1}. $$
Then
$$ \left(
E_H + z_0' V^+_{z_0';C(\lambda;z_0')}
\right)^{-1} = \frac{1}{z_0(|z_0'|^2-1)}
\left(
(z_0'-z_0) E_H \right.
$$
$$ \left.
+
z_0'(|z_0|^2-1)
\left(
E_H + z_0 V^+_{z_0;C(\lambda;z_0)}
\right)^{-1}
\right) $$
\begin{equation}
\label{f3_10}
=
\frac{1-z_0'\overline{z_0}}{1-|z_0'|^2}
\left(
E_H + \frac{z_0-z_0'}{1-z_0'\overline{z_0}} V^+_{z_0;C(\lambda;z_0)}
\right)
\left(
E_H + z_0 V^+_{z_0;C(\lambda;z_0)}
\right)^{-1},\ \lambda\in \mathbb{D}.
\end{equation}
\begin{lem}
\label{l3_1}
Let $z_0,z_0'\in \mathbb{D}$. Then
\begin{equation}
\label{f3_11}
\left|
\frac{z_0-z_0'}{1-z_0'\overline{z_0}}
\right|<1.
\end{equation}
\end{lem}
{\bf Proof. }
Consider the linear fractional transformation: $w=w(u)=\frac{z_0-u}{1-\overline{z_0}u}$.
If $|u|=1$ then $|1-\overline{z_0}u|=|u(\overline{u}-\overline{z_0})|=|u-z_0|$.
Moreover $w(z_0)=0$. Therefore $w$ maps $\mathbb{D}$ onto $\mathbb{D}$.
$\Box$

By~(\ref{f3_10}),(\ref{f3_11}) we may write:
$$ E_H + z_0' V^+_{z_0';C(\lambda;z_0')} $$
\begin{equation}
\label{f3_12}
=
\frac{1-|z_0'|^2}{1-z_0'\overline{z_0}}
\left( E_H + z_0 V^+_{z_0;C(\lambda;z_0)} \right)
\left(
E_H + \frac{z_0-z_0'}{1-z_0'\overline{z_0}} V^+_{z_0;C(\lambda;z_0)}
\right)^{-1},\quad \lambda\in \mathbb{D}.
\end{equation}
By~(\ref{f1_5}) we may write:
$$ V^+_{z_0;C(\lambda;z_0)} = V_{z_0} \oplus C(\lambda;z_0),\quad \lambda\in \mathbb{D}. $$
By conditions of the theorem it easily follows that $C(\widehat\lambda;z_0)$ is a contraction, and
\begin{equation}
\label{f3_13}
\exists u.-\lim_{k\to\infty} V^+_{z_0;C(\lambda_k;z_0)} = V_{z_0} \oplus C(\widehat\lambda;z_0)
= V^+_{z_0;C(\widehat\lambda;z_0)}.
\end{equation}
We may write
$$ \left\|
\left(
E_H + \frac{z_0-z_0'}{1-z_0'\overline{z_0}} V^+_{z_0;C(\lambda_k;z_0)}
\right)^{-1}
-
\left(
E_H + \frac{z_0-z_0'}{1-z_0'\overline{z_0}} V^+_{z_0;C(\widehat\lambda;z_0)}
\right)^{-1}
\right\|
$$
$$ \leq
\left\|
\left(
E_H + \frac{z_0-z_0'}{1-z_0'\overline{z_0}} V^+_{z_0;C(\lambda_k;z_0)}
\right)^{-1}
\right\| $$
\begin{equation}
\label{f3_14}
* \left|
\frac{z_0-z_0'}{1-z_0'\overline{z_0}} \right|
\left\|
V^+_{z_0;C(\lambda_k;z_0)} - V^+_{z_0;C(\widehat\lambda;z_0)}
\right\|
\left\|
\left(
E_H + \frac{z_0-z_0'}{1-z_0'\overline{z_0}} V^+_{z_0;C(\widehat\lambda;z_0)}
\right)^{-1}
\right\|.
\end{equation}
Since
$$ \left| \frac{z_0-z_0'}{1-z_0'\overline{z_0}} \right|
\left\| V^+_{z_0;C(\lambda_k;z_0)} \right\| \leq \delta < 1, $$
then
$$ \left\|
\left(
E_H + \frac{z_0-z_0'}{1-z_0'\overline{z_0}} V^+_{z_0;C(\lambda_k;z_0)}
\right) h
\right\|
\geq \left| \| h \| - \left\|
\frac{z_0-z_0'}{1-z_0'\overline{z_0}} V^+_{z_0;C(\lambda_k;z_0)}
\right|
\right| $$
$$ \geq (1-\delta)\| h \|; $$
$$ \left\|
\left(
E_H + \frac{z_0-z_0'}{1-z_0'\overline{z_0}} V^+_{z_0;C(\lambda_k;z_0)}
\right)^{-1} \right\| \leq \frac{1}{1-\delta}. $$
Passing to the limit in~(\ref{f3_14}) we obtain that
\begin{equation}
\label{f3_15}
u.-\lim_{k\to\infty} \left( E_H + \frac{z_0-z_0'}{1-z_0'\overline{z_0}} V^+_{z_0;C(\lambda_k;z_0)}
\right)^{-1} =
\left( E_H + \frac{z_0-z_0'}{1-z_0'\overline{z_0}} V^+_{z_0;C(\widehat\lambda;z_0)}
\right)^{-1}.
\end{equation}
By~(\ref{f3_15}),(\ref{f3_13}),(\ref{f3_12}) we conclude that there exists
\begin{equation}
\label{f3_16}
u.-\lim_{k\to\infty} V^+_{z_0';C(\lambda_k;z_0')}
= u.-\lim_{k\to\infty} V_{z_0'} \oplus C(\lambda_k;z_0')
=: V',
\end{equation}
such that
$$ E_H + z_0' V' $$
\begin{equation}
\label{f3_17}
=
\frac{1-|z_0'|^2}{1-z_0'\overline{z_0}}
\left( E_H + z_0 V^+_{z_0;C(\widehat\lambda;z_0)} \right)
\left(
E_H + \frac{z_0-z_0'}{1-z_0'\overline{z_0}} V^+_{z_0;C(\widehat\lambda;z_0)}
\right)^{-1}.
\end{equation}
By~(\ref{f3_16}) we see that $V'|_{M_{z_0'}} = V_{z_0'}$. Set
$$ C(\widehat\lambda;z_0') = V'|_{N_{z_0'}}. $$
Then relation~(\ref{f3_16}) shows that $C(\widehat\lambda;z_0')$ is a linear contraction which
maps $N_{z_0'}$ into $N_{\frac{1}{ \overline{z_0'} }}$. Thus, we have
\begin{equation}
\label{f3_18}
V' = V_{z_0'}\oplus C(\widehat\lambda;z_0') = V^+_{z_0';C(\widehat\lambda;z_0')}.
\end{equation}
By~(\ref{f3_16}),(\ref{f3_18}) we easily get that
\begin{equation}
\label{f3_18_1}
u.-\lim_{k\to\infty} C(\lambda_k;z_0')
= C(\widehat\lambda;z_0'),
\end{equation}
and~(\ref{f3_7}) is proved.

By~(\ref{f3_17}),(\ref{f3_18}) we obtain:
$$ \left( E_H + z_0' V^+_{z_0';C(\widehat\lambda;z_0')} \right)^{-1} $$
\begin{equation}
\label{f3_19}
= \frac{1-z_0'\overline{z_0}}{1-|z_0'|^2}
\left(
E_H + \frac{z_0-z_0'}{1-z_0'\overline{z_0}} V^+_{z_0;C(\widehat\lambda;z_0)}
\right)
\left( E_H + z_0 V^+_{z_0;C(\widehat\lambda;z_0)} \right)^{-1};
\end{equation}
$$ (1-|z_0'|^2) \left( E_H + z_0' V^+_{z_0';C(\widehat\lambda;z_0')} \right)^{-1} $$
$$ = \left(
(1-z_0'\overline{z_0}) E_H + (z_0-z_0') V^+_{z_0;C(\widehat\lambda;z_0)}
\right)
\left( E_H + z_0 V^+_{z_0;C(\widehat\lambda;z_0)} \right)^{-1};
$$
By subtracting $E_H$ from the both sides of  the last relation and by division by $-z_0'$
we get
$$ \frac{1}{z_0'}E_H + \frac{|z_0'|^2-1}{z_0'} \left( E_H + z_0' V^+_{z_0';C(\widehat\lambda;z_0')} \right)^{-1} $$
$$ = -\frac{1}{z_0'}\left(
(1-z_0'\overline{z_0}) E_H + (z_0-z_0') V^+_{z_0;C(\widehat\lambda;z_0)}
- ( E_H + z_0 V^+_{z_0;C(\widehat\lambda;z_0)} )
\right) $$
$$*
\left( E_H + z_0 V^+_{z_0;C(\widehat\lambda;z_0)} \right)^{-1}
$$
$$ = \left( \overline{z_0} E_H + V^+_{z_0;C(\widehat\lambda;z_0)} \right)
\left( E_H + z_0 V^+_{z_0;C(\widehat\lambda;z_0)} \right)^{-1}. $$
By~(\ref{f1_6}) and (\ref{f1_6_1}) we conclude that
\begin{equation}
\label{f3_20}
V_{C(\widehat\lambda;z_0');z_0'} = V_{C(\widehat\lambda;z_0);z_0},\qquad \forall z_0'\in \mathbb{D}.
\end{equation}
$\Box$

The following theorem is an analog of Theorem~2.1 (B) in~\cite{cit_5500_M}.
\begin{thm}
\label{t3_2}
Let $V$ be a closed isometric operator in a Hilbert space $H$, and $z_0\in \mathbb{D}\backslash\{ 0 \}$ be
a fixed point. Let $\mathbf{R}_z = \mathbf{R}_z(V)$ be an arbitrary
generalized resolvent of $V$, and $C(\lambda;z_0)\in \mathcal{S}(N_{z_0};N_{\frac{1}{ \overline{z_0} }})$
corresponds to $\mathbf{R}_z(V)$ by Inin's formula~(\ref{f1_7}).
$\mathbf{R}_z(V)$ has an analytic continuation to the set $\mathbb{D}\cup\mathbb{D}_e\cup \Delta$,
for some open arc $\Delta$ of $\mathbb{T}$, if and only if the following conditions are satisfied:

\begin{itemize}
\item[1)] $C(\lambda;z_0)$ has an extension to the set $\mathbb{D}\cup\Delta$ which is continuous in the
uniform operator topology;

\item[2)] The extended $C(\lambda;z_0)$ maps isometrically $N_{z_0}$ on the whole
$N_{ \frac{1}{ \overline{z_0} } }$, for all $\lambda\in\Delta$;

\item[3)] The operator $(E_H - \lambda V_{C(\lambda;z_0);z_0})^{-1}$ exists and it is defined on
the whole $H$, for all $\lambda\in\Delta$.
\end{itemize}
\end{thm}
{\bf Proof. }{\it Necessity.}
Choose an arbitrary point $\widehat\lambda\in\Delta$.
Let $z\in \mathbb{D}\backslash\{ 0 \}$ be an arbitrary point, and
$C(\lambda;z)\in \mathcal{S}(N_{z};N_{\frac{1}{ \overline{z} }})$
corresponds to the  generalized resolvent $\mathbf{R}_z(V)$ by Inin's formula~(\ref{f1_7}).
Using Inin's formula we may write:
$$ E_H - z V_{C(\lambda;z);z} = \frac{z}{\lambda} (E_H - \lambda V_{C(\lambda;z);z}) +
\left( 1-\frac{z}{\lambda}\right) E_H $$
$$ = \frac{z}{\lambda} \mathbf{R}_\lambda^{-1}
+ \left( 1-\frac{z}{\lambda}\right) E_H =
\frac{z}{\lambda} \left[
E_H + \left( \frac{\lambda}{z} - 1 \right) \mathbf{R}_\lambda
\right] \mathbf{R}_\lambda^{-1},\ \forall\lambda\in \mathbb{D}\backslash\{ 0 \}.
$$
Therefore
$$ \left[
E_H + \left( \frac{\lambda}{z} - 1 \right) \mathbf{R}_\lambda
\right]
=
\frac{\lambda}{z} ( E_H - z V_{C(\lambda;z);z} ) \mathbf{R}_\lambda, $$
has a bounded inverse defined on the whole $H$:
\begin{equation}
\label{f3_20_1}
\left[
E_H + \left( \frac{\lambda}{z} - 1 \right) \mathbf{R}_\lambda
\right]^{-1}
=
\frac{z}{\lambda} \mathbf{R}_\lambda^{-1} ( E_H - z V_{C(\lambda;z);z} )^{-1},\quad
\lambda \in \mathbb{D}\backslash\{ 0 \}.
\end{equation}
Then
\begin{equation}
\label{f3_21}
(E_H - z V_{C(\lambda;z);z})^{-1} = \frac{\lambda}{z}
\mathbf{R}_\lambda
\left[
E_H + \left( \frac{\lambda}{z} - 1 \right) \mathbf{R}_\lambda
\right]^{-1},\quad \lambda \in \mathbb{D}\backslash\{ 0 \}.
\end{equation}
Choose an arbitrary $\delta$: $0<\delta<1$. Assume that $z\in \mathbb{D}\backslash\{ 0 \}$
satisfies the following additional condition:
\begin{equation}
\label{f3_22}
\left| \frac{\widehat\lambda}{z} - 1 \right| \| \mathbf{R}_{\widehat\lambda} \| < \delta.
\end{equation}
Let us check that such points exist. If $\| \mathbf{R}_{\widehat\lambda} \|=0$, it is obvious.
In the opposite case, we look for $z=\varepsilon\widehat\lambda$,
$0<\varepsilon <1$. In this case
condition~(\ref{f3_22}) means that
$$ \left| \frac{1}{\varepsilon} - 1 \right| <
\frac{\delta}{ \| \mathbf{R}_{\widehat\lambda} \| }; $$
or $\varepsilon > \frac{1}{ 1 +
\frac{\delta}{ \| \mathbf{R}_{\widehat\lambda} \| } }$.

\noindent
Then there exists $\left[
E_H + \left( \frac{\widehat\lambda}{z} - 1 \right) \mathbf{R}_{\widehat\lambda}
\right]^{-1}$ which is bounded and defined on the whole $H$.
Moreover, by continuity inequality
\begin{equation}
\label{f3_22_1}
\left| \frac{\lambda}{z} - 1 \right| \| \mathbf{R}_{\lambda} \| < \delta,
\end{equation}
holds in an open neighborhood $U(\widehat\lambda)$
of $\widehat\lambda$, and therefore there exists
\begin{equation}
\label{f3_23}
\left[
E_H + \left( \frac{\lambda}{z} - 1 \right) \mathbf{R}_{\lambda}
\right]^{-1},\quad \forall\lambda\in U(\widehat\lambda),
\end{equation}
which is bounded and defined on the whole $H$.
We may write
$$ \left\|
\left( E_H + \left( \frac{\lambda}{z} - 1 \right) \mathbf{R}_{\lambda} \right) h
\right\| \geq
\left|
\| h \| - \left| \frac{\lambda}{z} - 1 \right|
\| \mathbf{R}_{\lambda} h \|
\right| $$
$$ \geq (1-\delta) \| h \|,\qquad h\in H. $$
Therefore
\begin{equation}
\label{f3_24}
\left\|
\left[
E_H + \left( \frac{\lambda}{z} - 1 \right) \mathbf{R}_{\lambda}
\right]^{-1}
\right\| \leq \frac{1}{1-\delta},\quad
\lambda\in U(\widehat\lambda).
\end{equation}
Choose an arbitrary sequence $\{ \lambda_k \}_{k=1}^\infty$
of points in $\mathbb{D}$, such that $\lambda_k\rightarrow\widehat\lambda$, as $k\rightarrow\infty$.
There exists a number $k_0\in \mathbb{N}$ such that $\lambda_k\in U(\widehat\lambda)\cap \mathbb{D}$,
$k\geq k_0$.
We may write
$$ \left\|
\left[
E_H + \left( \frac{\lambda_k}{z} - 1 \right) \mathbf{R}_{\lambda_k}
\right]^{-1}
-
\left[
E_H + \left( \frac{\widehat\lambda}{z} - 1 \right) \mathbf{R}_{\widehat\lambda}
\right]^{-1} \right\| $$
$$ \leq
\left \| \left[
E_H + \left( \frac{\lambda_k}{z} - 1 \right) \mathbf{R}_{\lambda_k}
\right]^{-1}
\right\|
\left\|
\left( \frac{\lambda_k}{z} - 1 \right) \mathbf{R}_{\lambda_k} -
\left( \frac{\widehat\lambda}{z} - 1 \right) \mathbf{R}_{\widehat\lambda}
\right\| $$
$$ *
\left\|
\left[
E_H + \left( \frac{\widehat\lambda}{z} - 1 \right) \mathbf{R}_{\widehat\lambda}
\right]^{-1}
\right\|. $$
The first factor on the right of the last equality is uniformly bounded by~(\ref{f3_24}).
Thus, we get
\begin{equation}
\label{f3_25}
u.-\lim_{\lambda\in \mathbb{D},\ \lambda\to\widehat\lambda}
\left[
E_H + \left( \frac{\lambda}{z} - 1 \right) \mathbf{R}_{\lambda}
\right]^{-1} =
\left[
E_H + \left( \frac{\widehat\lambda}{z} - 1 \right) \mathbf{R}_{\widehat\lambda}
\right]^{-1}.
\end{equation}
By relations~(\ref{f3_21}),(\ref{f3_25}) we conclude that
\begin{equation}
\label{f3_26}
u.-\lim_{\lambda\in \mathbb{D},\ \lambda\to\widehat\lambda}
(E_H - z V_{C(\lambda;z);z})^{-1} = \frac{\widehat\lambda}{z}
\mathbf{R}_{\widehat\lambda}
\left[
E_H + \left( \frac{\widehat\lambda}{z} - 1 \right) \mathbf{R}_{\widehat\lambda}
\right]^{-1}.
\end{equation}
Then there exists the following limit:
$$ u.-\lim_{\lambda\in \mathbb{D},\ \lambda\to\widehat\lambda}
V^+_{z;C(\lambda;z)} =
u.-\lim_{\lambda\in \mathbb{D},\ \lambda\to\widehat\lambda}
\left(
-\frac{1}{z} E_H + \frac{1-|z|^2}{z} (E_{H} - z V_{C(\lambda;z);z})^{-1}
\right) $$
\begin{equation}
\label{f3_27}
=
-\frac{1}{z} E_H +
\frac{1-|z|^2}{z^2}
\widehat\lambda
\mathbf{R}_{\widehat\lambda}
\left[
E_H + \left( \frac{\widehat\lambda}{z} - 1 \right) \mathbf{R}_{\widehat\lambda}
\right]^{-1}
=: V_z',
\end{equation}
where we used~(\ref{f1_6_2}). Notice that
$V^+_{z;C(\lambda;z)} = V_z\oplus C(\lambda;z)$.
Set
$C(\widehat\lambda;z) = V_z'|_{N_z}$. By~(\ref{f3_27}) we conclude that
$C(\widehat\lambda;z)$ is a linear contraction which maps
$N_z$ into $N_{\frac{1}{ \overline{z} }}$. Moreover, $V_z'|_{M_z} = V_z$, and therefore
$$ V_z' = V_z\oplus C(\widehat\lambda;z) = V^+_{z;C(\widehat\lambda;z)}. $$
By~(\ref{f3_27}) we easily obtain that
\begin{equation}
\label{f3_28}
u.-\lim_{\lambda\in \mathbb{D},\ \lambda\to\widehat\lambda}
C(\lambda;z) = C(\widehat\lambda;z).
\end{equation}
By Theorem~\ref{t3_1} we conclude that the last relation also holds for $z_0$, where
$C(\widehat\lambda;z_0)$ is a linear contraction
which maps $N_{z_0}$ into
$N_{\frac{1}{ \overline{z_0} }}$, and
$V_{C(\widehat\lambda;z_0);z_0} = V_{C(\widehat\lambda;z);z}$.

\noindent
Comparing relation~(\ref{f3_27}) for $V_z' = V^+_{z;C(\widehat\lambda;z)}$, with
relation~(\ref{f1_6_2}) we get:
\begin{equation}
\label{f3_28_1}
\frac{\widehat\lambda}{z}
\mathbf{R}_{\widehat\lambda}
\left[
E_H + \left( \frac{\widehat\lambda}{z} - 1 \right) \mathbf{R}_{\widehat\lambda}
\right]^{-1} = \left( E_H - z V_{C(\widehat\lambda;z);z} \right)^{-1},
\end{equation}
for the above choice of $z$.

Thus, we have extended by continuity the function $C(\lambda;z_0)$ to the set
$\mathbb{D}\cup\Delta$. Let us check that this extension is continuous in the uniform
operator topology. It remains to check that for an arbitrary $\widehat\lambda\in\Delta$ we have
\begin{equation}
\label{f3_29}
u.-\lim_{\lambda\in \mathbb{D}\cup\Delta,\ \lambda\to\widehat\lambda}
C(\lambda;z_0) = C(\widehat\lambda;z_0).
\end{equation}
We choose $z\in \mathbb{D}\backslash\{ 0 \}$ satisfying~(\ref{f3_22}) and construct
a neighborhood $U(\widehat\lambda)$, as before.
For an arbitrary $\lambda\in U(\widehat\lambda)$ we may write:
$$ \left\|
\left[
E_H + \left( \frac{\lambda}{z} - 1 \right) \mathbf{R}_{\lambda}
\right]^{-1}
-
\left[
E_H + \left( \frac{\widehat\lambda}{z} - 1 \right) \mathbf{R}_{\widehat\lambda}
\right]^{-1}
\right\| $$
$$ \leq
\left\|
\left[
E_H + \left( \frac{\lambda}{z} - 1 \right) \mathbf{R}_{\lambda}
\right]^{-1}
\right\|
\left\|
\left( \frac{\widehat\lambda}{z} - 1 \right) \mathbf{R}_{\widehat\lambda} -
\left( \frac{\lambda}{z} - 1 \right) \mathbf{R}_{\lambda}
\right\| $$
$$ *
\left\|
\left[
E_H + \left( \frac{\widehat\lambda}{z} - 1 \right) \mathbf{R}_{\widehat\lambda}
\right]^{-1}
\right\|. $$
By~(\ref{f3_24}) we obtain that
\begin{equation}
\label{f3_30}
u.-\lim_{\lambda\to\widehat\lambda}
\left[
E_H + \left( \frac{\lambda}{z} - 1 \right) \mathbf{R}_{\lambda}
\right]^{-1}
=
\left[
E_H + \left( \frac{\widehat\lambda}{z} - 1 \right) \mathbf{R}_{\widehat\lambda}
\right]^{-1}.
\end{equation}
By~(\ref{f3_28_1}) we conclude that
\begin{equation}
\label{f3_31}
\frac{\lambda}{z}
\mathbf{R}_{\lambda}
\left[
E_H + \left( \frac{\lambda}{z} - 1 \right) \mathbf{R}_{\lambda}
\right]^{-1} = \left( E_H - z V_{C(\lambda;z);z} \right)^{-1},\quad \forall\lambda\in
(U(\widehat\lambda)\cap\overline{\mathbb{D}})\backslash\{ 0 \},
\end{equation}
for the above choice of $z$. In fact, for an arbitrary $\widetilde\lambda\in\Delta\cap U(\widehat\lambda)$,
there exists a neighborhood $\widetilde U(\widetilde\lambda)\subset U(\widehat\lambda)$, where
inequality~(\ref{f3_22_1}) holds for the same choice of $z$.
Then repeating the arguments after~(\ref{f3_22_1}) for $\widetilde\lambda$ instead of $\widehat\lambda$,
we obtain that~(\ref{f3_31}) holds for $\widetilde\lambda$.
For the points inside $\mathbb{D}$ we may use relation~(\ref{f3_21}).

By relations~(\ref{f3_30}),(\ref{f3_31}) we get
\begin{equation}
\label{f3_32}
u.-\lim_{\lambda\in \mathbb{D}\cup\Delta,\ \lambda\to\widehat\lambda}
\left( E_H - z V_{C(\lambda;z);z} \right)^{-1}
=
\left( E_H - z V_{C(\widehat\lambda;z);z} \right)^{-1}.
\end{equation}
Since it was proven that $V_{C(\widehat\lambda;z);z}$ does not depend on the choice of
$z\in \mathbb{D}\backslash\{ 0 \}$ (and for $z\in \mathbb{D}$ this fact follows from Inin's formula),
the last relation holds for all $z\in \mathbb{D}\backslash\{ 0 \}$.

By relation~(\ref{f1_6_2}) we obtain that
\begin{equation}
\label{f3_33}
u.-\lim_{\lambda\in \mathbb{D}\cup\Delta,\ \lambda\to\widehat\lambda}
V^+_{z;C(\lambda;z)}
= V^+_{z;C(\widehat\lambda;z)},\quad \forall z\in \mathbb{D}\backslash\{ 0 \},
\end{equation}
and therefore relation~(\ref{f3_29}) holds.
Thus, condition~1) in the statement of the theorem is proven.

By~(\ref{f3_28_1}) we see that
\begin{equation}
\label{f3_34}
\mathbf{R}_{\widehat\lambda}
=
\frac{z}{\widehat\lambda}
\left( E_H - z V_{C(\widehat\lambda;z);z} \right)^{-1}
\left[
E_H + \left( \frac{\widehat\lambda}{z} - 1 \right) \mathbf{R}_{\widehat\lambda}
\right],
\end{equation}
for the above choice of $z$. Therefore $\mathbf{R}_{\widehat\lambda}$ has a bounded inverse, defined on
the whole $H$. Then
$$ E_H = \frac{z}{\widehat\lambda}
\left( E_H - z V_{C(\widehat\lambda;z);z} \right)^{-1}
\left[
E_H + \left( \frac{\widehat\lambda}{z} - 1 \right) \mathbf{R}_{\widehat\lambda}
\right]
\mathbf{R}_{\widehat\lambda}^{-1}; $$
$$ E_H - z V_{C(\widehat\lambda;z);z} =
\frac{z}{\widehat\lambda}
\left[
E_H + \left( \frac{\widehat\lambda}{z} - 1 \right) \mathbf{R}_{\widehat\lambda}
\right]
\mathbf{R}_{\widehat\lambda}^{-1} $$
$$ = \frac{z}{\widehat\lambda}
\mathbf{R}_{\widehat\lambda}^{-1} + \left( 1 - \frac{z}{\widehat\lambda} \right) E_H. $$
From the last relation we get
\begin{equation}
\label{f3_35}
\mathbf{R}_{\widehat\lambda} =
\left( E_H - \widehat\lambda V_{C(\widehat\lambda;z);z} \right)^{-1}.
\end{equation}
Since $V_{C(\widehat\lambda;z);z}$ does not depend on the choice of $z$,
the last relation holds for all $z\in \mathbb{D}\backslash\{ 0 \}$.
Consequently, condition~3) in the statement of the theorem is proven.

By the property~(\ref{f3_3}), passing to the limit as $z\rightarrow\widehat\lambda$, we get
\begin{equation}
\label{f3_36}
\mathbf{R}_{\widehat\lambda}^* = E_H - \mathbf{R}_{\widehat\lambda}.
\end{equation}
On the other hand, by~(\ref{f3_35}) we get
\begin{equation}
\label{f3_37}
\mathbf{R}_{\widehat\lambda}^* = \left( E_H - \overline{\widehat\lambda}
V_{C(\widehat\lambda;z);z}^* \right)^{-1}.
\end{equation}
By~(\ref{f3_35})-(\ref{f3_37}) we see that
\begin{equation}
\label{f3_38}
E_H = \left( E_H - \overline{\widehat\lambda}
V_{C(\widehat\lambda;z);z}^* \right)^{-1}
+
\left( E_H - \widehat\lambda V_{C(\widehat\lambda;z);z} \right)^{-1},\quad
\forall z\in \mathbb{D}\backslash\{ 0 \}.
\end{equation}
By multiplying the both sides of the last relation by $(E_H - \overline{\widehat\lambda}
V_{C(\widehat\lambda;z);z}^*)$ from the left, and by
$(E_H - \widehat\lambda V_{C(\widehat\lambda;z);z})$ from the right, we get
$$ (E_H - \overline{\widehat\lambda}
V_{C(\widehat\lambda;z);z}^*)(E_H - \widehat\lambda V_{C(\widehat\lambda;z);z}) =
E_H - \widehat\lambda V_{C(\widehat\lambda;z);z}
+
E_H - \overline{\widehat\lambda}
V_{C(\widehat\lambda;z);z}^*. $$
After multiplication in the left-hand side and simplification we obtain that
$$ V_{C(\widehat\lambda;z);z}^* V_{C(\widehat\lambda;z);z} = E_H,\qquad \forall z\in \mathbb{D}\backslash\{ 0 \}. $$
On the other hand, by~(\ref{f3_38}) we may write:
$$ \left( E_H - \overline{\widehat\lambda}
V_{C(\widehat\lambda;z);z}^* \right)^{-1}
= E_H - \left( E_H - \widehat\lambda V_{C(\widehat\lambda;z);z} \right)^{-1} $$
$$ = -\widehat\lambda V_{C(\widehat\lambda;z);z}
\left( E_H - \widehat\lambda V_{C(\widehat\lambda;z);z} \right)^{-1}; $$
$$ V_{C(\widehat\lambda;z);z} = -\frac{1}{ \widehat\lambda }
\left( E_H - \overline{\widehat\lambda}
V_{C(\widehat\lambda;z);z}^* \right)^{-1}
\left( E_H - \widehat\lambda V_{C(\widehat\lambda;z);z} \right). $$
Since $\left( E_H - \widehat\lambda V_{C(\widehat\lambda;z);z} \right)^{-1}$ is defined on the whole $H$
and bounded, we conclude that $R(V_{C(\widehat\lambda;z);z}) = H$.
Hence, the operator $V_{C(\widehat\lambda;z);z}$ is unitary in $H$.
Then the corresponding operator $V^+_{z;C(\widehat\lambda;z)} = V_z\oplus
C(\widehat\lambda;z)$ is unitary, as well.
In particular, this fact implies that $C(\widehat\lambda;z)$ is isometric and maps
$N_z$ on the whole $N_{ \frac{1}{\overline{z}} }$. Since $z$ is an arbitrary point from
$\mathbb{D}\backslash\{ 0 \}$, we obtain that condition~2) of the theorem is satisfied.

{\it Sufficiency. }
Let conditions 1)-3) be satisfied.
Choose an arbitrary $\widehat\lambda\in\Delta$.
Choose an arbitrary sequence $\{ \lambda_k \}_{k=1}^\infty$
of points in $\mathbb{D}\cup\Delta$, such that $\lambda_k\rightarrow\widehat\lambda$, as $k\rightarrow\infty$.
Using~(\ref{f1_6_1}) we write:
$$ E_H - \lambda_k V_{C(\lambda_k;z_0);z_0} =
\left( 1-\frac{1}{z_0} \right) E_H - \frac{|z_0|^2-1}{z_0} (E_{H} + z_0 V^+_{z_0;C(\lambda_k;z_0)})^{-1} $$
$$ = (1-\lambda_k \overline{z_0})
\left[
E_H + \frac{z_0-\lambda_k}{1-\lambda_k \overline{z_0}} V^+_{z_0;C(\lambda_k;z_0)}
\right]
(E_{H} + z_0 V^+_{z_0;C(\lambda_k;z_0)})^{-1}; $$
$$ E_H - \widehat\lambda V_{C(\widehat\lambda;z_0);z_0} = $$
\begin{equation}
\label{f3_39}
= (1-\widehat\lambda \overline{z_0})
\left[
E_H + \frac{z_0-\widehat\lambda}{1-\widehat\lambda \overline{z_0}} V^+_{z_0;C(\widehat\lambda;z_0)}
\right]
(E_{H} + z_0 V^+_{z_0;C(\widehat\lambda;z_0)})^{-1}.
\end{equation}
By~(\ref{f3_39}) we write:
$$
E_H + \frac{z_0-\widehat\lambda}{1-\widehat\lambda \overline{z_0}} V^+_{z_0;C(\widehat\lambda;z_0)}
=
\frac{1}{1-\widehat\lambda \overline{z_0}}
\left( E_H - \widehat\lambda V_{C(\widehat\lambda;z_0);z_0} \right)
\left( E_{H} + z_0 V^+_{z_0;C(\widehat\lambda;z_0)} \right). $$
Since $V_{C(\widehat\lambda;z_0);z_0}$ is closed, by condition~3) it follows that
there exists
$( E_H - \widehat\lambda V_{C(\widehat\lambda;z_0);z_0} )^{-1}$, which is defined on the whole $H$, and bounded.
Therefore there exists
$[ E_H + \frac{z_0-\widehat\lambda}{1-\widehat\lambda \overline{z_0}} V^+_{z_0;C(\widehat\lambda;z_0)} ]^{-1}$,
which is bounded and defined on the whole $H$.
From~(\ref{f3_39}) it follows that
$$ \left( E_H - \widehat\lambda V_{C(\widehat\lambda;z_0);z_0} \right)^{-1} $$
\begin{equation}
\label{f3_41}
=
\frac{1}{1-\widehat\lambda \overline{z_0}}
(E_{H} + z_0 V^+_{z_0;C(\widehat\lambda;z_0)})
\left[
E_H + \frac{z_0-\widehat\lambda}{1-\widehat\lambda \overline{z_0}} V^+_{z_0;C(\widehat\lambda;z_0)}
\right]^{-1}.
\end{equation}
For points $\lambda_k$ which belong to $\Delta$ we may apply the same argument, while for
points $\lambda_k$ from $\mathbb{D}$ we can use Lemma~\ref{l3_1}, to obtain an analogous
representation:
$$ \left( E_H - \lambda_k V_{C(\lambda_k;z_0);z_0} \right)^{-1} $$
\begin{equation}
\label{f3_41_1}
= \frac{1}{1-\lambda_k \overline{z_0}}
(E_{H} + z_0 V^+_{z_0;C(\lambda_k;z_0)})
\left[
E_H + \frac{z_0-\lambda_k}{1-\lambda_k \overline{z_0}} V^+_{z_0;C(\lambda_k;z_0)}
\right]^{-1},\quad k\in \mathbb{N}.
\end{equation}

Observe that by condition~1) we have:
$$ \left\|
V^+_{z_0;C(\lambda_k;z_0)} - V^+_{z_0;C(\widehat\lambda;z_0)}
\right\|
=
\sup_{h\in H,\ \| h \| = 1} \left\|
(C(\lambda_k;z_0) - C(\widehat\lambda;z_0)) P_{N_{z_0}} h
\right\| $$
$$ \leq
\left\|
C(\lambda_k;z_0) - C(\widehat\lambda;z_0)
\right\| \rightarrow 0,\quad k\rightarrow\infty; $$
\begin{equation}
\label{f3_42}
u.-\lim_{k\rightarrow\infty} V^+_{z_0;C(\lambda_k;z_0)} = V^+_{z_0;C(\widehat\lambda;z_0)}.
\end{equation}
Let us check that there exists an open neighborhood $U_1(\widehat\lambda)$ of $\widehat\lambda$, and a number
$K>0$ such that
\begin{equation}
\label{f3_43}
\left\|
\left[
E_H + \frac{z_0-\lambda_k}{1-\lambda_k \overline{z_0}} V^+_{z_0;C(\lambda_k;z_0)}
\right]^{-1}
\right\| \leq K,\quad
\forall\lambda_k:\ \lambda_k\in U_1(\widehat\lambda).
\end{equation}
The latter condition is equivalent to the following condition:
\begin{equation}
\label{f3_44}
\left\|
\left[
E_H + \frac{z_0-\lambda_k}{1-\lambda_k \overline{z_0}} V^+_{z_0;C(\lambda_k;z_0)}
\right] g
\right\| \geq \frac{1}{K} \| g \|,\quad \forall g\in H,\
\forall\lambda_k:\ \lambda_k\in U_1(\widehat\lambda).
\end{equation}
Suppose to the contrary that condition~(\ref{f3_44}) is not true.
Choose a sequence of open balls $U^n(\widehat\lambda)$ with the centrum at $\widehat\lambda$
and radius $\frac{1}{n}$; and set $K_n = n$, $n\in \mathbb{N}$.
Then for each $n\in \mathbb{N}$, there exist elements $g_n\in H$, and
$\lambda_{k_n}\in U^n(\widehat\lambda)$,
$k_n\in \mathbb{N}$, such that:
\begin{equation}
\label{f3_45}
\left\|
\left[
E_H + \frac{z_0-\lambda_{k_n}}{1-\lambda_{k_n} \overline{z_0}} V^+_{z_0;C(\lambda_{k_n};z_0)}
\right] g_n
\right\| < \frac{1}{n} \| g_n \|.
\end{equation}
It is clear that $g_n$ are all non-zero. Set $\widehat g_n = \frac{g_n}{\| g_n \|_H}$, $n\in \mathbb{N}$.
Then
\begin{equation}
\label{f3_46}
\left\|
\left[
E_H + \frac{z_0-\lambda_{k_n}}{1-\lambda_{k_n} \overline{z_0}} V^+_{z_0;C(\lambda_{k_n};z_0)}
\right] \widehat g_n
\right\| < \frac{1}{n}.
\end{equation}
Since $|\lambda_{k_n} - \widehat\lambda| < \frac{1}{n}$, then
$\lim_{n\to\infty} \lambda_{k_n} = \widehat\lambda$.
Then we may write
$$ \frac{1}{n} >
\left\|
\left[
E_H + \frac{z_0-\widehat\lambda}{1-\widehat\lambda \overline{z_0}} V^+_{z_0;C(\widehat\lambda;z_0)}
\right] \widehat g_n \right. $$
$$\left. +
\left(
\frac{z_0-\lambda_{k_n}}{1-\lambda_{k_n} \overline{z_0}} V^+_{z_0;C(\lambda_{k_n};z_0)}
-
\frac{z_0-\widehat\lambda}{1-\widehat\lambda \overline{z_0}} V^+_{z_0;C(\widehat\lambda;z_0)}
\right)
\widehat g_n
\right\| $$
$$ \geq \left|
\left\|
\left[
E_H + \frac{z_0-\widehat\lambda}{1-\widehat\lambda \overline{z_0}} V^+_{z_0;C(\widehat\lambda;z_0)}
\right] \widehat g_n \right\| \right. $$
$$ \left. -
\left\|
\frac{z_0-\lambda_{k_n}}{1-\lambda_{k_n} \overline{z_0}} V^+_{z_0;C(\lambda_{k_n};z_0)}
-
\frac{z_0-\widehat\lambda}{1-\widehat\lambda \overline{z_0}} V^+_{z_0;C(\widehat\lambda;z_0)}
\right\|
\right|
$$
\begin{equation}
\label{f3_47}
\geq
L
-
\left\|
\frac{z_0-\lambda_{k_n}}{1-\lambda_{k_n} \overline{z_0}} V^+_{z_0;C(\lambda_{k_n};z_0)}
-
\frac{z_0-\widehat\lambda}{1-\widehat\lambda \overline{z_0}} V^+_{z_0;C(\widehat\lambda;z_0)}
\right\|,\ L>0,
\end{equation}
for sufficiently large $n$,
since $[E_H + \frac{z_0-\widehat\lambda}{1-\widehat\lambda \overline{z_0}} V^+_{z_0;C(\widehat\lambda;z_0)}]$
has a bounded inverse on the whole $H$, and
the norm in the right-hand side tends to zero. Passing to the limit in relation~(\ref{f3_47}) as
$n\rightarrow\infty$, we obtain a contradiction.

\noindent
Consequently, there exists an open neighborhood $U_1(\widehat\lambda)$ of $\widehat\lambda$, and a number
$K>0$ such that inequality~(\ref{f3_43}) holds.
We may write:
$$ \left\|
\left[
E_H + \frac{z_0-\lambda_k}{1-\lambda_k \overline{z_0}} V^+_{z_0;C(\lambda_k;z_0)}
\right]^{-1}
-
\left[
E_H + \frac{z_0-\widehat\lambda}{1-\widehat\lambda \overline{z_0}} V^+_{z_0;C(\widehat\lambda;z_0)}
\right]^{-1}
\right\| $$
$$ \leq
\left\|
\left[
E_H + \frac{z_0-\lambda_k}{1-\lambda_k \overline{z_0}} V^+_{z_0;C(\lambda_k;z_0)}
\right]^{-1}
\right\|
\left\|
\frac{z_0-\widehat\lambda}{1-\widehat\lambda \overline{z_0}} V^+_{z_0;C(\widehat\lambda;z_0)}
-
\frac{z_0-\lambda_k}{1-\lambda_k \overline{z_0}} V^+_{z_0;C(\lambda_k;z_0)}
\right\|
$$
$$ * \left\|
\left[
E_H + \frac{z_0-\widehat\lambda}{1-\widehat\lambda \overline{z_0}} V^+_{z_0;C(\widehat\lambda;z_0)}
\right]^{-1}
\right\|. $$
By~(\ref{f3_43}) we conclude that
\begin{equation}
\label{f3_48}
u.-\lim_{k\rightarrow\infty} \left[
E_H + \frac{z_0-\lambda_k}{1-\lambda_k \overline{z_0}} V^+_{z_0;C(\lambda_k;z_0)}
\right]^{-1}
=
\left[
E_H + \frac{z_0-\widehat\lambda}{1-\widehat\lambda \overline{z_0}} V^+_{z_0;C(\widehat\lambda;z_0)}
\right]^{-1}.
\end{equation}
By~(\ref{f3_41}),(\ref{f3_41_1}),(\ref{f3_42}),(\ref{f3_48}) we conclude that
\begin{equation}
\label{f3_49}
u.-\lim_{k\rightarrow\infty}
\left( E_H - \lambda_k V_{C(\lambda_k;z_0);z_0} \right)^{-1}
=
\left( E_H - \widehat\lambda V_{C(\widehat\lambda;z_0);z_0} \right)^{-1},
\end{equation}
and therefore
\begin{equation}
\label{f3_50}
u.-\lim_{\lambda\in \mathbb{D}\cup\Delta,\ \lambda\to\widehat\lambda}
\left( E_H - \lambda V_{C(\lambda;z_0);z_0} \right)^{-1}
=
\left( E_H - \widehat\lambda V_{C(\widehat\lambda;z_0);z_0} \right)^{-1}.
\end{equation}
By Inin's formula, for $\lambda\in \mathbb{D}$, we have
$\left( E_H - \lambda V_{C(\lambda;z_0);z_0} \right)^{-1} = \mathbf{R}_\lambda$.
Thus, relation~(\ref{f3_50}) shows that the operator-valued function $\mathbf{R}_\lambda$,
$\lambda\in \mathbb{D}$, has a continuation to the set $\mathbb{D}\cup\Delta$, which is
continuous in the uniform operator topology.

Choose an arbitrary $h\in H$ and consider an analytic function
\begin{equation}
\label{f3_50_1}
f(\lambda) = f_h(\lambda) = (\mathbf{R}_\lambda h,h),\qquad \lambda\in \mathbb{D}.
\end{equation}
Then $f(\lambda)$ has a continuous extension to $\mathbb{D}\cup \Delta$, which is equal to
$$ f(\lambda) = \left( \left( E_H - \lambda V_{C(\lambda;z_0);z_0} \right)^{-1} h,h \right),\quad
\lambda\in \mathbb{D}\cup \Delta. $$
Let us check that
\begin{equation}
\label{f3_51}
\left( E_H - \lambda V_{C(\lambda;z_0);z_0} \right)^{-1}
= E_H -
\left( E_H - \overline{\lambda} V^*_{C(\lambda;z_0);z_0} \right)^{-1},\quad
\forall\lambda\in\Delta.
\end{equation}
Choose an arbitrary $\lambda\in\Delta$.
By condition~2) of the theorem we conclude that $V_{C(\lambda;z_0);z_0}$ is unitary.
Then
\begin{equation}
\label{f3_52}
(E_H - \overline{\lambda}
V_{C(\lambda;z);z}^*)(E_H - \lambda V_{C(\lambda;z);z}) =
E_H - \lambda V_{C(\lambda;z);z}
+
E_H - \overline{\lambda}
V_{C(\lambda;z);z}^*.
\end{equation}
To verify the last relation, it is sufficient to make multiplication in the left-hand side and
a simplification. Multiplying~(\ref{f3_52}) by $\left( E_H - \overline{\lambda} V^*_{C(\lambda;z_0);z_0} \right)^{-1}$
from the left, and by $\left( E_H - \lambda V_{C(\lambda;z_0);z_0} \right)^{-1}$ from the right,
we easily get~(\ref{f3_51}).

We may write:
$$ \overline{f(\lambda)} = \left( h,\left( E_H - \lambda V_{C(\lambda;z_0);z_0} \right)^{-1} h \right) =
\left( \left( E_H - \overline{\lambda} V^*_{C(\lambda;z_0);z_0} \right)^{-1} h,h \right) $$
$$ = (h,h) - \left(\left( E_H - \lambda V_{C(\lambda;z_0);z_0} \right)^{-1} h,h \right) =
(h,h) - f(\lambda);$$
\begin{equation}
\label{f3_53}
{\mathrm Re} f(\lambda) = \frac{1}{2} (h,h),\qquad  \lambda\in \Delta.
\end{equation}
Set $g(\lambda) = if(\lambda) - \frac{i}{2} (h,h)$, $\lambda\in \mathbb{D}\cup\Delta$.
Then ${\mathrm Im} g(\lambda)=0$.
Consequently, by the Schwarc principle, $g(\lambda)$ admits an analytic continuation
$\widetilde g(\lambda) = \widetilde g_h(\lambda)$
to the set $\mathbb{D}\cup\Delta\cup \mathbb{D}_e$.
Moreover, we have
\begin{equation}
\label{f3_54}
\widetilde g(\lambda) = \overline{ \widetilde g\left( \frac{1}{ \overline{\lambda} } \right) },\qquad
\lambda\in \mathbb{D}_e.
\end{equation}
Then
$$ \widetilde f(\lambda) = \widetilde f_h(\lambda) := \frac{1}{i} \widetilde g(\lambda) + \frac{1}{2}(h,h),\quad
\lambda\in \mathbb{D}\cup\Delta\cup \mathbb{D}_e, $$
is an analytic continuation of $f(\lambda)$. By~(\ref{f3_54}) we get
\begin{equation}
\label{f3_55}
\widetilde f(\lambda) = - \overline{ \widetilde f\left( \frac{1}{ \overline{\lambda} } \right) }
+ (h,h),\qquad
\lambda\in \mathbb{D}_e.
\end{equation}
Using~(\ref{f3_3}) we may write:
$$ \widetilde f(\lambda) = - \overline{ \left( \mathbf{R}_{\frac{1}{ \overline{\lambda} }} h,h \right)_H }
+ (h,h)_H =
- \left( h, ( E_H - \mathbf{R}^*_{\lambda} ) h \right)_H
+ (h,h)_H
$$
\begin{equation}
\label{f3_56}
= (\mathbf{R}_\lambda h,h)_H,\qquad \lambda\in \mathbb{D}_e.
\end{equation}
Set
$$ R_\lambda (h,g) = \frac{1}{4} \left( \widetilde f_{h+g}(\lambda) - \widetilde f_{h-g}(\lambda) +
i \widetilde f_{h+ig}(\lambda) - i \widetilde f_{h-ig}(\lambda) \right), $$
\begin{equation}
\label{f3_57}
h,g\in H,\quad \lambda\in \mathbb{D}\cup\Delta\cup \mathbb{D}_e.
\end{equation}
Observe that $R_\lambda (h,g)$ is an analytic function of $\lambda$ on $\mathbb{D}\cup\Delta\cup \mathbb{D}_e$.
By~(\ref{f3_50_1}),(\ref{f3_56}) we see that
\begin{equation}
\label{f3_58}
R_\lambda (h,g) = (\mathbf{R}_\lambda h,g)_H,\quad
h,g\in H,\quad \lambda\in \mathbb{D}\cup \mathbb{D}_e.
\end{equation}
By~(\ref{f3_50}) we see that
$$ R_{\widehat\lambda} (h,g) =
\lim_{\lambda\in \mathbb{D},\ \lambda\to\widehat\lambda}
(\mathbf{R}_\lambda h,g)_H $$
$$ =
\lim_{\lambda\in \mathbb{D},\ \lambda\to\widehat\lambda}
\left( \left( E_H - \lambda V_{C(\lambda;z_0);z_0} \right)^{-1} h, g\right)_H $$
\begin{equation}
\label{f3_59}
=
\left( \left( E_H - \widehat\lambda V_{C(\widehat\lambda;z_0);z_0} \right)^{-1} h, g\right)_H,\quad
h,g\in H,\quad \widehat\lambda\in\Delta.
\end{equation}
Then the following operator-valued function:
$$ T_\lambda = \left\{
\begin{array}{cc} \mathbf{R}_\lambda, & \lambda\in \mathbb{D}\cup \mathbb{D}_e \\
\left( E_H - \lambda V_{C(\lambda;z_0);z_0} \right)^{-1} & \lambda\in\Delta
\end{array}\right., $$
is an extension of $\mathbf{R}_\lambda$, which is analytic in the weak operator topology,
and therefore in the uniform operator topology.
$\Box$

\begin{cor}
\label{c3_1}
Theorem~\ref{t3_2} remains valid for the choice $z_0=0$.
\end{cor}
{\bf Proof. }
Let $V$ be a closed isometric operator in a Hilbert space $H$.
Let $\mathbf{R}_z(V)$ be an arbitrary
generalized resolvent of $V$, and $F(\lambda) = C(\lambda;0)\in
\mathcal{S}(N_{0};N_{\infty})$
corresponds to $\mathbf{R}_z(V)$ by Inin's formula~(\ref{f1_7}) for $z_0=0$, which
in this case coincides with Chumakin's formula~(\ref{f1_1}).
Consider an arbitrary open arc $\Delta$ of $\mathbb{T}$.

Choose an arbitrary point $z_0\in \mathbb{D}\backslash\{ 0 \}$.
Consider the following isometric operator
$$ \mathbf{V} = (V+\overline{z_0}E_H)( E_H + z_0 V )^{-1},\quad D(\mathbf{V})=(E_H+z_0 V)D(V). $$
Then
$$ V = (\mathbf{V}-\overline{z_0}E_H)(E_H - z_0 \mathbf{V})^{-1} = \mathbf{V}_{z_0}. $$
Recall that generalized resolvents of $\mathbf{V}$ and $\mathbf{V}_{z_0}$ are related
by~(\ref{f2_2}) and this correspondence is bijective.
Let $\mathbf{R}_z(\mathbf{V})$ be the generalized resolvent which corresponds by~(\ref{f2_2}) to
the generalized resolvent $\mathbf{R}_z(\mathbf{V}_{z_0})=\mathbf{R}_z(V)$.

By~(\ref{f2_2}) we see that $\mathbf{R}_{\widetilde t}(\mathbf{V}_{z_0})$ has a limit as
$\widetilde t\rightarrow \widetilde{t}_0\in\Delta$, if and only if
$\mathbf{R}_{\widetilde u}(\mathbf{V})$ has a limit as
$\widetilde u\rightarrow \widetilde{u}_0\in\Delta_1$, where
$$ \Delta_1 = \left\{ \widetilde u:\ \widetilde u = \frac{\widetilde t + z_0}
{1+\overline{z_0} \widetilde t},\ \widetilde t\in\Delta \right\}. $$
Thus, $\mathbf{R}_{\widetilde t}(\mathbf{V}_{z_0})$ has an extension by continuity to
$\mathbb{T}_e\cup\Delta$, iff $\mathbf{R}_{\widetilde u}(\mathbf{V})$
has an extension by continuity to
$\mathbb{T}_e\cup\Delta_1$. The extended values are related by~(\ref{f2_2}), as well.
From~(\ref{f2_2}) we see that the extension of $\mathbf{R}_{\widetilde t}(\mathbf{V}_{z_0})$ is
analytic iff the extension of $\mathbf{R}_{\widetilde u}(\mathbf{V})$ is analytic.
Consequently, $\mathbf{R}_{\widetilde t}(V) = \mathbf{R}_{\widetilde t}(\mathbf{V}_{z_0})$ has an analytic extension
to
$\mathbb{T}_e\cup\Delta$, if and only if  $\mathbf{R}_{\widetilde u}(\mathbf{V})$
has an analytic extension to $\mathbb{T}_e\cup\Delta_1$.

By~Theorem~\ref{t3_2}, $\mathbf{R}_{\widetilde u}(\mathbf{V})$
has an analytic extension to $\mathbb{T}_e\cup\Delta_1$ iff

\begin{itemize}
\item[1)] $C(\lambda;z_0)$ has an extension to the set $\mathbb{D}\cup\Delta_1$ which is continuous in the
uniform operator topology;

\item[2)] The extended $C(\lambda;z_0)$ maps isometrically $N_{z_0}$ on the whole
$N_{ \frac{1}{ \overline{z_0} } }$, for all $\lambda\in\Delta_1$;

\item[3)] The operator $(E_H - \lambda \mathbf{V}_{C(\lambda;z_0);z_0})^{-1}$ exists and it is defined on
the whole $H$, for all $\lambda\in\Delta_1$,
\end{itemize}
where $C(\lambda;z_0)\in \mathcal{S}(N_{z_0};N_{\frac{1}{ \overline{z_0} }})$
corresponds to $\mathbf{R}_z(\mathbf{V})$ by Inin's formula~(\ref{f1_7}).
Recall that $C(\lambda;z_0)$ is related to $F(\widetilde t)$ in the following way:
$$ C(\widetilde u) = F\left( \frac{\widetilde u - z_0}{1-\overline{z_0} \widetilde u} \right),\quad
u\in \mathbb{T}_e. $$
By using this relation we easily get that condition~1) is equivalent to

1') $F(\widetilde t)$ has an extension to the set $\mathbb{D}\cup\Delta$ which is continuous in the
uniform operator topology;

and condition~2) is equivalent to

2') The extended $F(\widetilde t)$ maps isometrically $N_0$ on the whole
$N_{\infty}$, for all $\widetilde t \in\Delta$.

By~(\ref{f3_39}) we conclude that
$(E_H - \lambda \mathbf{V}_{C(\lambda;z_0);z_0})^{-1}$ exists and is defined on
the whole $H$,  for all $\lambda\in\Delta_1$, iff
$$ \left[ E_H - \frac{\lambda - z_0}{1-\lambda\overline{z_0}} \left( \mathbf{V}_{z_0}\oplus C(\lambda;z_0) \right)
\right]^{-1} =
\left[ E_H - \widetilde t \left( \mathbf{V}_{z_0}\oplus F(\widetilde t) \right)
\right]^{-1}, $$
exists and is defined on the whole $H$,  for all $\widetilde t\in\Delta$.
$\Box$

\section{Some decompositions of a Hilbert space.}
The following result appeared in~\cite[Lemma 1]{cit_5950_R}. However, its proof was based on
Shmulyan's lemma. It seems that no correct proof of this lemma ever appeared published.
An attempt to prove Shmulyan's lemma was performed by L.A.~Shtraus in~\cite[Lemma]{cit_5900_S_LA}.
However, the proof was not complete. We shall use the idea of L.A.~Shtraus to prove the
following weaker result.
\begin{thm}
\label{t4_1}
Let $V$ be a closed isometric operator in a Hilbert space $H$.
Let $\zeta\in \mathbb{T}$,  and $\zeta^{-1}$ be a point of  the regular type of $V$.
Then the following decompositions are valid:
\begin{equation}
\label{f4_1}
D(V) \dotplus N_\zeta = H;
\end{equation}
\begin{equation}
\label{f4_2}
R(V) \dotplus N_\zeta = H.
\end{equation}
\end{thm}
{\bf Proof. }
At first, we suppose that $\zeta = 1$.
Let us check that
\begin{equation}
\label{f4_3}
\| Vf + g \|_H = \| f+g \|_H,\qquad f\in D(V),\ g\in N_1.
\end{equation}
In fact, we may write:
$$ \| Vf + g \|_H^2 = (Vf,Vf)_H + (Vf,g)_H + (g,Vf)_H + (g,g)_H $$
$$ = \| f \|_H^2 + (Vf,g)_H + (g,Vf)_H + \|g\|_H^2. $$
Since $0 = ((E-V)f, g)_H = (f,g)_H - (Vf,g)_H$, we get
\begin{equation}
\label{f4_4}
(Vf,g)=(f,g),\qquad f\in D(V),\ g\in N_1,
\end{equation}
and therefore
$$ \| Vf + g \|_H^2 =
\| f \|_H^2 + (f,g)_H + (g,f)_H + \|g\|_H^2 = \| f+g \|_H^2. $$
Consider the following operator:
\begin{equation}
\label{f4_5}
U(f+g) = Vf+g,\qquad f\in D(V),\ g\in N_1.
\end{equation}
Let us check that this operator is well-defined, with the domain $D(U)= D(V) + N_1$.
Let $h\in D(U)$ has two representations:
$$ h = f_1+g_1 = f_2+g_2,\qquad f_1,f_2\in D(V),\ g_1,g_2\in N_1. $$
Then using~(\ref{f4_3}) we may write
$$ \| Vf_1+g_1 - (Vf_2+g_2) \|_H^2 =  \| V(f_1-f_2) + (g_1-g_2) \|_H^2 =
\| f_1-f_2 + g_1-g_2 \|_H^2 = 0. $$
Thus, $U$ is well-defined. Moreover, $U$ is linear and using~(\ref{f4_4}) we write:
$$ (U(f+g),U(\widetilde f + \widetilde g)) =
(Vf+g,V\widetilde f + \widetilde g) = (Vf,V\widetilde f) +
(Vf,\widetilde g) + (g,V\widetilde f) + (g,\widetilde g) $$
$$ = (f,\widetilde f) + (f,\widetilde g) + (g,\widetilde f) + (g,\widetilde g) =
(f+g,\widetilde f + \widetilde g), $$
for $f,\widetilde f\in D(V)$, $g,\widetilde g\in N_1$.
Therefore $U$ is isometric.

Suppose that there exists $h\in H$, $h\not= 0$, $h\in D(V)\cap N_1$. Then
$$ 0 = U 0 = U( h + (-h)) = Vh - h = (V-E_H)h, $$
and this contradicts to the fact that $1$ is a point of the regular type of $V$.
Therefore
\begin{equation}
\label{f4_6}
D(V)\cap N_1 = \{ 0 \}.
\end{equation}
Notice that we do not know, a priori, that $D(U)$ is closed.
Consider the following operator $W$:
$$ W = U|_S, $$
where
$$ S = \left\{ h\in D(U):\ h\perp N_1 \right\} = D(U)\cap M_1. $$
Thus, $W$ is an isometric operator with the domain $D(W) = D(U)\cap M_1$.
Choose an arbitrary element $g\in D(U)$. Then
$g = g_{M_1} + g_{N_1}$, $g_{M_1}\in M_1$, $g_{N_1}\in N_1\subseteq D(U)$.
Therefore $g_{M_1} = P^H_{M_1} g \in D(U)$, $g_{M_1}\perp N_1$. Thus, $g_{M_1}\in D(W)$;
\begin{equation}
\label{f4_7}
P^H_{M_1} D(U) \subseteq D(W).
\end{equation}
On the other hand, choose an arbitrary $h\in D(W)$. Then $h\in D(U)\cap M_1$, and therefore
$h = P^H_{M_1} h\in P^H_{M_1} D(U)$. Consequently, we have
\begin{equation}
\label{f4_8}
D(W) = P^H_{M_1} D(U) = P^H_{M_1} (D(V)+N_1) = P^H_{M_1}D(V) \subseteq M_1.
\end{equation}
Choose an arbitrary $f\in D(V)$. Let $f= f_{M_1} + f_{N_1}$, $f_{M_1}\in M_1$, $f_{N_1}\in N_1$.
Then $f-f_{N_1}\in D(U)$, and $f-f_{N_1}\perp N_1$, i.e. $f-f_{N_1}\in D(W)$.
We may write
$$ (W-E_H)(f-f_{N_1}) = (U-E_H)(f-f_{N_1}) = U(f-f_{N_1}) - f + f_{N_1} = Vf - f; $$
\begin{equation}
\label{f4_9}
(W-E_H) D(W) \supseteq (V-E_H)D(V) = M_1.
\end{equation}
On the other hand, choose an arbitrary $w\in D(W)$, $w = w_{D(V)} + w_{N_1}$, $w_{D(V)}\in D(V)$,
$w_{N_1}\in N_1$.
Since $w\perp N_1$, we may write:
\begin{equation}
\label{f4_9_1}
w = P^H_{M_1} w = P^H_{M_1} w_{D(V)}.
\end{equation}
We may write
$$ (W-E_H) w = Uw - w = Vw_{D(V)} + w_{N_1} -w_{D(V)} - w_{N_1} = (V-E_H) w_{D(V)}; $$
$$ (W-E_H) D(W) \subseteq (V-E_H) D(V) = M_1. $$
From the last relation and~(\ref{f4_9}) we obtain:
\begin{equation}
\label{f4_10}
(W-E_H) D(W) = (V-E_H)D(V) = M_1.
\end{equation}
Moreover, if $(W-E_H) w = 0$, then $(V-E_H) w_{D(V)}=0$; and therefore $w_{D(V)}=0$. By~(\ref{f4_9_1})
this implies $w=0$. Consequently, there exists $(W-E_H)^{-1}$. By~(\ref{f4_10}) we get
\begin{equation}
\label{f4_11}
D(W) = (W-E_H)^{-1} M_1.
\end{equation}
Since $D(W)\subseteq M_1$, using~(\ref{f4_10}), we may write:
$$ W w = (W-E_H)w + w\in M_1; $$
\begin{equation}
\label{f4_12}
D(W)\subseteq M_1,\quad W D(W) \subseteq M_1.
\end{equation}
Consider the closure $\overline{W}$ of $W$ with the domain $\overline{D(W)}$.
By~(\ref{f4_12}) we see that
\begin{equation}
\label{f4_13}
D(\overline{W})\subseteq M_1,\quad \overline{W} D(\overline{W}) \subseteq M_1.
\end{equation}
Then
$$ (\overline{W} - E_H) D(\overline{W}) \subseteq M_1. $$
On the other hand, by~(\ref{f4_10}) we have
$$ (\overline{W}  - E_H) D(\overline{W}) \supseteq (\overline{W}  - E_H) D(W) =
(W  - E_H) D(W) = M_1. $$
Thus, we conclude that
\begin{equation}
\label{f4_14}
(\overline{W} - E_H) D(\overline{W}) = M_1.
\end{equation}
Let us check that there exists $(\overline{W} - E_H)^{-1}$. Suppose to the contrary that
there exists $h\in D(\overline{W})$, $h\not=0$, such that $(\overline{W} - E_H)h = 0$.
By the definition of the closure, there exists a sequence $h_n\in D(W)$, $n\in \mathbb{N}$,
which converges to $h$, and
$W h_n\rightarrow \overline{W}h$, as $n\rightarrow\infty$.
Then $(W-E_H)h_n \rightarrow \overline{W}h - h = 0$, as $n\rightarrow\infty$.
Let $h_n = h_{1;n} + h_{2;n}$, $h_{1;n}\in D(V)$, $h_{2;n}\in N_1$, $n\in \mathbb{N}$.
Then $(W-E_H)h_n = U h_n - h_n = Vh_{1;n} + h_{2;n} - h_{1;n} - h_{2;n} = (V-E_H) h_{1;n}$.
Therefore $(V-E_H) h_{1;n} \rightarrow 0$, as $n\rightarrow\infty$.
Since $(V-E_H)$ has a bounded inverse, we conclude that $h_{1;n} \rightarrow 0$, as $n\rightarrow\infty$.
Then $h_n = P^H_{M_1} h_n = P^H_{M_1} h_{1;n} \rightarrow 0$, as $n\rightarrow\infty$.
Therefore, we get $h=0$, and this contradicts to our assumption.

Thus, there exists $(\overline{W} - E_H)^{-1}\supseteq (W - E_H)^{-1}$.
By~(\ref{f4_14}),(\ref{f4_10}), we conclude that
$$ (\overline{W} - E_H)^{-1} = (W - E_H)^{-1}, $$
and therefore
$$ \overline{W} = W. $$
Thus $W$ may be viewed as a closed isometric operator in a Hilbert space $M_1$.
Then $(W - E_H)^{-1}$ is closed, and it is defined on $M_1$.
Therefore $(W - E_H)^{-1}$ is bounded. This means that $1$ is a regular point of $W$.
Therefore $W$ is a unitary operator in $M_1$.
In particular, this implies that $D(W)=R(W)=M_1$.

By the definition of $W$ we conclude that $D(W) = M_1\subseteq D(U)$, and $U|_{M_1} = W$.
On the other hand, $U|_{N_1} = E_{N_1}$. Then $D(U)=H$, and
$$ U = W\oplus E_{N_1}. $$
Thus, $U$ is a unitary operator. In particular, $D(U)=R(U)=H$, which means that
\begin{equation}
\label{f4_15}
D(V)+N_1=H,\quad R(V)+N_1=H.
\end{equation}
The first sum is direct by~(\ref{f4_6}).
Suppose that $h\in R(V)\cap N_1$.  Then $h = Vf$, $f\in D(V)$, and we may write:
$$ 0 = Vf + (-h) = U(f+(-h)). $$
Since $U$ is unitary, we get $f=h=Vf$, $(V-E_H)f=0$, and therefore $f=0$, and $h=0$.
Thus, the second sum in~(\ref{f4_15}) is direct, as well.
So, we proved the theorem for the case $\zeta = 1$.

In the general case, we can apply the proven part of the theorem to
$\widehat V := \zeta V$.
$\Box$

\begin{cor}
\label{c4_1}
In conditions of Theorem~\ref{t4_1} the following decompositions are valid:
\begin{equation}
\label{f4_16}
\overline{ (H\ominus D(V)) \dotplus M_\zeta } = H;
\end{equation}
\begin{equation}
\label{f4_17}
\overline{ (H\ominus R(V)) \dotplus M_\zeta } = H.
\end{equation}
\end{cor}
{\bf Proof. }
The proof is based on the following lemma.
\begin{lem}
\label{l4_1}
Let $M_1$ and $M_2$ be two subspaces in a Hilbert space $H$, such that
$M_1\cap M_2 = \{ 0 \}$, and
\begin{equation}
\label{f4_18}
M_1 \dotplus M_2 = H.
\end{equation}
Then
\begin{equation}
\label{f4_19}
\overline{ (H\ominus M_1) \dotplus (H\ominus M_2) } = H.
\end{equation}
\end{lem}
{\bf Proof. }
Suppose that $h\in H$, is such that $h\in ((H\ominus M_1) \cap (H\ominus M_2))$. Then
$h\perp M_1$, $h\perp M_2$, and therefore $h\perp (M_1+M_2)$, $h\perp H$, $h=0$.

Suppose that $g\in H$, is such that $g\perp ((H\ominus M_1) \dotplus (H\ominus M_2))$. Then
$g\in M_1$, $g\in M_2$, and therefore $g=0$.
$\Box$

By applying the lemma with $M_1 = D(V)$, $M_2 = N_\zeta$, and $M_1 = R(V)$, $M_2 = N_\zeta$,
we complete the proof of the corollary.
$\Box$

\section{Isometric operators with gaps.}
Now we shall present full proofs, some slight correction and generalization of Ryabtseva results~\cite{cit_5950_R}.

For the sake of convenience we put here the proofs of the following lemmas, see Lemma 2 and its corollary, and
Lemma~3 in~\cite{cit_5950_R}.
\begin{lem}
\label{l5_1}
Let $V$ be a closed isometric operator in a Hilbert space $H$, and $\zeta\in \mathbb{T}$.
Then the following equality holds:
\begin{equation}
\label{f5_1}
V P^H_{M_0} f = \zeta^{-1} P^H_{M_\infty} f,\qquad \forall f\in N_\zeta,
\end{equation}
and therefore
\begin{equation}
\label{f5_2}
\left\| P^H_{M_0} f \right\| =
\left\| P^H_{M_\infty} f \right\|,\qquad \forall f\in N_\zeta;
\end{equation}
\begin{equation}
\label{f5_3}
\left\| P^H_{N_0} f \right\| =
\left\| P^H_{N_\infty} f \right\|,\qquad \forall f\in N_\zeta.
\end{equation}
\end{lem}
{\bf Proof. } Choose an arbitrary $f\in N_\zeta$.
For an arbitrary $u\in D(V)=M_0$ we may write:
$$ \left( \zeta^{-1} f - V P^H_{D(V)} f, Vu \right)_H
= \zeta^{-1} (f,Vu)_H - \left( P^H_{D(V)} f, u \right)_H $$
$$ = (f,\zeta V u)_H - (f,u)_H = \left( f, (\zeta V - E_H) u \right)_H = 0. $$
Therefore $( \zeta^{-1} f - V P^H_{M_0} f )\perp M_\infty$.
By applying $P^H_{M_\infty}$ to this element we get~(\ref{f5_1}).
Relation~(\ref{f5_2}) is obvious, since $V$ is isometric, and~(\ref{f5_3}) easily follows.
$\Box$

\begin{lem}
\label{l5_2}
Let $V$ be a closed isometric operator in a Hilbert space $H$, and
$C$ be a linear bounded operator in $H$, $D(C) = N_0$, $R(C)\subseteq N_\infty$.
Let $\zeta\in \mathbb{T}$, and $\zeta^{-1}$ be an eigenvalue of $V^+_{0;C} = V\oplus C$.
If $f\in H$, $f\not= 0$, is an eigenvector of $V^+_{0;C}$ corresponding to $\zeta^{-1}$, then
$f\in N_\zeta$, and
\begin{equation}
\label{f5_4}
C P^H_{N_0} f = \zeta^{-1} P^H_{N_\infty} f.
\end{equation}
\end{lem}
{\bf Proof. }
Let $f$ be an eigenvector of $V^+_{0;C}$ corresponding to $\zeta^{-1}\in \mathbb{T}$:
$$ (V\oplus C) f = V P^H_{M_0} f + C P^H_{N_0} f = \zeta^{-1} \left( P^H_{M_\infty} f + P^H_{N_\infty} f
\right). $$
By the orthogonality of the summands we see that the last relation is equivalent to the following
relations
\begin{equation}
\label{f5_5}
V P^H_{M_0} f = \zeta^{-1} P^H_{M_\infty} f;
\end{equation}
\begin{equation}
\label{f5_6}
C P^H_{N_0} f = \zeta^{-1} P^H_{N_\infty} f,
\end{equation}
and therefore~(\ref{f5_4}) follows.
Relation~(\ref{f5_5}) implies that
$P^H_{M_\infty} (\zeta^{-1} f - V P^H_{M_0} f) = 0$;
$( \zeta^{-1} f - V P^H_{M_0} f )\perp M_\infty$. Then for arbitrary $u\in D(V)$, we may write:
$$ 0 = \left( \zeta^{-1} f - V P^H_{M_0} f , Vu \right)_H
= \zeta^{-1} (f,Vu)_H -
\left( P^H_{M_0} f , u \right)_H $$
$$ = (f,\zeta Vu)_H - (f,u)_H
= (f, (\zeta V-E_H)u )_H. $$
Therefore $f\in N_\zeta$.
$\Box$

For an arbitrary $\zeta\in \mathbb{T}$, we define an operator $W_\zeta$ by the following equality:
\begin{equation}
\label{f5_7}
W_\zeta P^H_{N_0} f = \zeta^{-1} P^H_{N_\infty} f,\qquad f\in N_\zeta,
\end{equation}
with the domain $D(W_\zeta) = P^H_{N_0} N_\zeta$.
Let us check that this definition is correct. If $g\in D(W_\zeta)$ admits two representations:
$g = P^H_{N_0} f_1 = P^H_{N_0} f_2$, $f_1,f_2\in N_\zeta$, then
$P^H_{N_0} (f_1-f_2) = 0$. By~(\ref{f5_3}) this implies
$P^H_{N_\infty} (f_1-f_2) = 0$, and therefore the definition is correct.
The operator $W_\zeta$ is linear, and
$$ \left\| W_\zeta P^H_{N_0} f \right\|_H =
\left\| P^H_{N_\infty} f \right\|_H =
\left\| P^H_{N_0} f \right\|_H. $$
Thus, $W_\zeta$ is isometric. Notice that $R(W_\zeta) = P^H_{N_\infty} N_\zeta$.

\noindent
Set
$$ S = P^H_{N_0}|_{N_\zeta},\quad Q = P^H_{N_\infty}|_{N_\zeta}. $$
In what follows, we suppose that {\it $\zeta^{-1}$ is a point of the regular type of the operator $V$}.
Let us check that in this case operators $S$ and $Q$ are invertible.
Suppose to the contrary that there exists $f\in N_\zeta$, $f\not=0$: $Sf = P^H_{N_0} f = 0$.
Then $f = P^H_{M_0} f \not= 0$. Thus, we get by Theorem~\ref{t4_1} that
$f\in M_0\cap N_\zeta = \{ 0 \}$. We obtained a contradiction.
In a similar way, suppose that there exists $g\in N_\zeta$, $g\not=0$: $Qf = P^H_{N_\infty} f = 0$.
Then $g = P^H_{M_\infty} g \not= 0$. Therefore by Theorem~\ref{t4_1} we get
$g\in M_\infty\cap N_\zeta = \{ 0 \}$. We obtained a contradiction, as well.

\noindent
Moreover, by Theorem~\ref{t4_1} we conclude that
\begin{equation}
\label{f5_8}
P^H_{N_0} N_\zeta = N_0;\quad P^H_{N_\infty} N_\zeta = N_\infty.
\end{equation}
Thus, $S^{-1}$ and $Q^{-1}$ are closed and defined on the subspaces $N_0$ and $N_\infty$, respectively.
Therefore $S^{-1}$ and $Q^{-1}$ are bounded.
By~(\ref{f5_7}) we see that $D(W_\zeta) = N_0$, $R(W_\zeta) = N_\infty$, and
\begin{equation}
\label{f5_8_1}
W_\zeta = \zeta^{-1} Q S^{-1}.
\end{equation}

The following theorem holds, see~\cite[Theorem 1]{cit_5950_R}.
\begin{thm}
\label{t5_1}
Let $V$ be a closed isometric operator in a Hilbert space $H$, and
$C$ be a linear bounded operator in $H$, $D(C) = N_0$, $R(C)\subseteq N_\infty$.
Let $\zeta\in \mathbb{T}$, and $\zeta^{-1}$ be a point of the regular type of the operator $V$.
The point $\zeta^{-1}$ is an eigenvalue of $V^+_{0;C} = V\oplus C$, if and only if
\begin{equation}
\label{f5_9}
(C - W_\zeta) g = 0,\qquad g\in N_0,\ g\not=0.
\end{equation}
\end{thm}
{\bf Proof. }
{\it Necessity. }
Since $\zeta^{-1}$ is an eigenvalue of $V^+_{0;C} = V\oplus C$, then by Lemma~\ref{l5_2} there exists
$f\in N_\zeta$, $f\not= 0$, such that
\begin{equation}
\label{f5_10}
C P^H_{N_0} f = \zeta^{-1} P^H_{N_\infty} f.
\end{equation}
Comparing the last relation with the definition of $W_\zeta$, we see that $C P^H_{N_0} f =
W_\zeta P^H_{N_0} f$. Set $g=P^H_{N_0} f = Sf$. Since $S$ is invertible, then $g\not= 0$.

{\it Sufficiency. } From~(\ref{f5_9}) we get~(\ref{f5_10}) with $f:= S^{-1}g$.
By~Lemma~\ref{l5_1} we see that relations~(\ref{f5_5}),(\ref{f5_6}) hold. The latter, as we have seen
before relation~(\ref{f5_5}), is equivalent to the fact that
$\zeta^{-1}$ is an eigenvalue of $V^+_{0;C} = V\oplus C$, with the eigenvector $f$.
$\Box$

The following theorem is a slightly corrected version of~\cite[Theorem 2]{cit_5950_R}.
\begin{thm}
\label{t5_2}
Let $V$ be a closed isometric operator in a Hilbert space $H$, and
$C$ be a linear bounded operator in $H$, $D(C) = N_0$, $R(C)\subseteq N_\infty$.
Let $\zeta\in \mathbb{T}$, and $\zeta^{-1}$ be a point of the regular type of the operator $V$.
Then
\begin{equation}
\label{f5_11}
R\left(
V^+_{0;C} - \zeta^{-1} E_H
\right) = H,
\end{equation}
if and only if the following two conditions hold:
\begin{equation}
\label{f5_12}
\left(
C - W_\zeta
\right) N_0 = N_\infty;
\end{equation}
\begin{equation}
\label{f5_12_1}
P^H_{M_\infty} M_\zeta = M_\infty.
\end{equation}
\end{thm}
{\bf Proof. }
{\it Necessity. }
Choose an arbitrary $h\in N_\infty$. By~(\ref{f5_11}) there exists $x\in H$, such that
\begin{equation}
\label{f5_14}
\left( V^+_{0;C} - \zeta^{-1} E_H \right) x = (V\oplus C) x - \zeta^{-1} x = h.
\end{equation}
For an arbitrary $u\in D(V)$, we may write:
$$ (x,(E_H-\zeta V)u)_H = (x,(V^{-1} - \zeta E_H) Vu)_H =
(x,((V^+_{0;C})^* - \zeta E_H) Vu)_H $$
$$ = ((V^+_{0;C} - \zeta^{-1}E_H)x,Vu)_H = (h,Vu)_H = 0, $$
and therefore $x\in N_\zeta$. Set $g=Sx\in N_0$, and using~(\ref{f5_8_1}) write:
$$ (C-W_\zeta) g = CSx - W_\zeta Sx = CSx - \zeta^{-1} Q x. $$
Since $h\in N_\infty$, we apply $P^H_{N_\infty}$ to the equality~(\ref{f5_14}) to get
$$ C P^H_{N_0} x - \zeta^{-1} P^H_{N_\infty} x = h; $$
$$ C S x - \zeta^{-1} Q x = h. $$
Therefore we obtain that
$$ (C-W_\zeta) g = h, $$
and~(\ref{f5_12}) holds.

Choose an arbitrary $\widehat h\in M_\infty$. By~(\ref{f5_11}) there exists $\widehat x\in H$, such that
\begin{equation}
\label{f5_15}
\left( V^+_{0;C} - \zeta^{-1} E_H \right) \widehat x = (V\oplus C) \widehat x - \zeta^{-1} \widehat x = \widehat h.
\end{equation}
The last equality is equivalent to the following two equalities obtained by applying projectors
$P^H_{M_\infty}$ and $P^H_{N_\infty}$:
\begin{equation}
\label{f5_16}
V P^H_{M_0} \widehat x - \zeta^{-1} P^H_{M_\infty} \widehat x = \widehat h;
\end{equation}
\begin{equation}
\label{f5_17}
C P^H_{N_0} \widehat x - \zeta^{-1} P^H_{N_\infty} \widehat x = 0.
\end{equation}
By Theorem~\ref{t4_1} we may write:
$$ \widehat x = x_{D(V)} + x_{N_\zeta},\qquad x_{D(V)}\in D(V),\ x_{N_\zeta}\in N_\zeta. $$
By substitution this decomposition into relation~(\ref{f5_17}) we get:
$$ C P^H_{N_0} x_{N_\zeta} - \zeta^{-1} P^H_{N_\infty} x_{N_\zeta} - \zeta^{-1} P^H_{N_\infty} x_{D(V)} = 0; $$
\begin{equation}
\label{f5_18}
(C - W_\zeta) P^H_{N_0} x_{N_\zeta} = \zeta^{-1} P^H_{N_\infty} x_{D(V)}.
\end{equation}
On the other hand, by substitution of the decomposition into~(\ref{f5_16}) we get:
$$ V x_{D(V)} + V P^H_{M_0} x_{N_\zeta} - \zeta^{-1} P^H_{M_\infty} x_{D(V)} - \zeta^{-1}
P^H_{M_\infty} x_{N_\zeta} = \widehat h; $$
\begin{equation}
\label{f5_19}
V x_{D(V)}  - \zeta^{-1} P^H_{M_\infty} x_{D(V)} = \widehat h,
\end{equation}
where we used~Lemma~\ref{l5_1}. Then
$$ P^H_{M_\infty} ( V - \zeta^{-1} E_H ) x_{D(V)} = \widehat h, $$
and~(\ref{f5_12_1}) follows directly.

{\it Sufficiency. }
Choose an arbitrary $h\in H$, $h = h_1 + h_2$, $h_1\in M_\infty$, $h_2\in N_\infty$. By~(\ref{f5_12}) there
exists $g\in N_0$ such that
$$ \left( C - W_\zeta \right) g = Cg - W_\zeta g = h_2. $$
Set $x = S^{-1} g\in N_\zeta$. Then
$$ CS x - \zeta^{-1} Q x = h_2; $$
$$ P^H_{N_\infty} (V\oplus C) P^H_{N_0} x  - \zeta^{-1} P^H_{N_\infty} x = h_2; $$
$$ P^H_{N_\infty} \left( V^+_{0;C} x - \zeta^{-1} x \right) = h_2. $$
Observe that by Lemma~\ref{l5_1} we may write:
$$ P^H_{M_\infty} \left( V^+_{0;C} x - \zeta^{-1} x \right) =
V P^H_{M_0} x - \zeta^{-1} P^H_{M_\infty} x = 0. $$
Therefore
\begin{equation}
\label{f5_20}
\left( V^+_{0;C} - \zeta^{-1} E_H \right) x = h_2.
\end{equation}
By~(\ref{f5_12_1}) there exists $w\in M_\zeta$, such that
$$ P^H_{M_\infty} w = h_1. $$
Let
$$ w = ( V - \zeta^{-1} E_H ) \widetilde x_{D(V)},\quad \widetilde x_{D(V)}\in D(V). $$
Then
\begin{equation}
\label{f5_21}
V \widetilde x_{D(V)}  - \zeta^{-1} P^H_{M_\infty} \widetilde x_{D(V)} = h_1.
\end{equation}
By~(\ref{f5_12}) there exists $r\in N_0$ such that
$$ \left( C - W_\zeta \right) r = \zeta^{-1} P^H_{N_\infty} \widetilde x_{D(V)}. $$
Set
$\widetilde x_{N_\zeta} := S^{-1} r\in N_\zeta$. Then
$$ (C-W_\zeta) P^H_{N_0} \widetilde x_{N_\zeta} =
\zeta^{-1} P^H_{N_\infty} \widetilde x_{D(V)}; $$
\begin{equation}
\label{f5_22}
C P^H_{N_0} \widetilde x_{N_\zeta} - \zeta^{-1} P^H_{N_\infty} \widetilde x_{N_\zeta}
- \zeta^{-1} P^H_{N_\infty} \widetilde x_{D(V)} = 0.
\end{equation}
Set $\widetilde x = \widetilde x_{D(V)} + \widetilde x_{N_\zeta}$.
Then from~(\ref{f5_22}) we get
\begin{equation}
\label{f5_23}
C P^H_{N_0} \widetilde x - \zeta^{-1} P^H_{N_\infty} \widetilde x = 0.
\end{equation}
Using~(\ref{f5_21}) and Lemma~\ref{l5_1} we write:
$$ h_1 = V \widetilde x_{D(V)}  - \zeta^{-1} P^H_{M_\infty} \widetilde x_{D(V)} +
V P^H_{M_0} \widetilde x_{N_\zeta} - \zeta^{-1} P^H_{M_\infty} \widetilde x_{N_\zeta} $$
\begin{equation}
\label{f5_24}
= V P^H_{M_0} \widetilde x - \zeta^{-1} P^H_{M_\infty} \widetilde x.
\end{equation}
Summing relations~(\ref{f5_23}) and~(\ref{f5_24}) we get
\begin{equation}
\label{f5_25}
(V\oplus C) \widetilde x - \zeta^{-1} \widetilde  x
= \left( V^+_{0;C} - \zeta^{-1} E_H \right) \widetilde  x
= h_1.
\end{equation}
Summing relations~(\ref{f5_20}) and~(\ref{f5_25}) we deduce that relation~(\ref{f5_11}) holds.
$\Box$

The following theorem holds, cf.~\cite[Theorem 4]{cit_5950_R}:
\begin{thm}
\label{t5_3}
Let $V$ be a closed isometric operator in a Hilbert space $H$, and
$\Delta$ be some open arc of $\mathbb{T}$, such that $\zeta^{-1}$ is a point of the regular type of $V$,
$\forall\zeta\in\Delta$.
Let the following condition be satisfied:
\begin{equation}
\label{f5_26}
P^H_{M_\infty} M_\zeta = M_\infty,\qquad \forall\zeta\in\Delta.
\end{equation}
Let $\mathbf{R}_z = \mathbf{R}_z(V)$ be an arbitrary
generalized resolvent of $V$, and $C(\lambda;0)\in \mathcal{S}(N_{0};N_\infty)$
corresponds to $\mathbf{R}_z(V)$ by Inin's formula~(\ref{f1_7}).
$\mathbf{R}_z(V)$ has an analytic continuation to the set $\mathbb{D}\cup\mathbb{D}_e\cup \Delta$,
if and only if the following conditions are satisfied:

\begin{itemize}
\item[1)] $C(\lambda;0)$ has an extension to the set $\mathbb{D}\cup\Delta$ which is continuous in the
uniform operator topology;

\item[2)] The extended $C(\lambda;0)$ maps isometrically $N_{0}$ on the whole
$N_\infty$, for all $\lambda\in\Delta$;

\item[3)] The operator $C(\lambda;0)-W_\lambda$ is invertible for all $\lambda\in\Delta$, and
\begin{equation}
\label{f5_27}
(C(\lambda;0)-W_\lambda) N_0 = N_\infty,\qquad \forall\lambda\in\Delta.
\end{equation}
\end{itemize}
\end{thm}
{\bf Proof. }{\it Necessity.}
Let $\mathbf{R}_z(V)$ have an analytic continuation to the set $\mathbb{D}\cup\mathbb{D}_e\cup \Delta$.
By Corollary~\ref{c3_1} we conclude that conditions~1) and 2) are satisfied, and
the operator $(E_H - \lambda V_{C(\lambda;0);0})^{-1} = -\frac{1}{\lambda}
( (V\oplus C(\lambda;0)) - \frac{1}{\lambda} E_H )^{-1}$ exists and it is defined on
the whole $H$, for all $\lambda\in\Delta$.
By Theorem~\ref{t5_1} we conclude that the operator
$C(\lambda;0)-W_\lambda$ is invertible for all $\lambda\in\Delta$.
By Theorem~\ref{t5_2} we obtain that relation~(\ref{f5_27}) holds.

{\it Sufficiency. } Let conditions 1)-3) be satisfied.
By Theorem~\ref{t5_2} we get that $R( (V\oplus C(\lambda;0)) - \frac{1}{\lambda} E_H ) =
R(E_H - \lambda V_{C(\lambda;0);0}) = H$.
By Theorem~\ref{t5_1} we see that the operator $(V\oplus C(\lambda;0)) - \frac{1}{\lambda} E_H$
is invertible. By Corollary~\ref{c3_1} we conclude that $\mathbf{R}_z(V)$ has
an analytic continuation to the set $\mathbb{D}\cup\mathbb{D}_e\cup \Delta$.
$\Box$

\begin{rmr}
\label{r5_1}
By Corollary~\ref{c3_1}, if $\mathbf{R}_z(V)$ has an analytic
continuation to the set $\mathbb{D}\cup\mathbb{D}_e\cup \Delta$, then
$( (V\oplus C(\lambda;0)) - \frac{1}{\lambda} E_H )^{-1}$ exists and it is bounded.
Therefore points $\lambda^{-1}$, $\lambda\in\Delta$, are of the regular type for $V$.
On the other hand, by Theorem~\ref{t5_2}, condition~(\ref{f5_26}) holds in this case.
Thus, by Proposition~\ref{p3_1}, condition~(\ref{f5_26}) and the condition that
points $\lambda^{-1}$, $\lambda\in\Delta$, are of the regular type for $V$,
both are {\it necessary} for the existence of a spectral function $\mathbf{F}$ of $V$ such that
$\mathbf{F}(\overline{\Delta}) = 0$. Thus, they do not imply on the generality of Theorem~\ref{t5_3}.
\end{rmr}

We shall obtain an analogous result in terms of the parameter $C(\lambda;z_0)$ of Inin's formula,
for an arbitrary $z_0\in \mathbb{D}$.
\begin{thm}
\label{t5_4}
Let $V$ be a closed isometric operator in a Hilbert space $H$, and
$\Delta$ be some open arc of $\mathbb{T}$, such that $\zeta^{-1}$ is a point of the regular type of $V$,
$\forall\zeta\in\Delta$. Let $z_0\in \mathbb{D}$ be an arbitrary fixed point, and
the following condition be satisfied:
\begin{equation}
\label{f5_28}
P^H_{M_{ \frac{1}{ \overline{z_0} } }} M_\frac{z_0 + \zeta}{1+\zeta \overline{z_0}} =
M_{ \frac{1}{ \overline{z_0} } } ,\qquad \forall\zeta\in\Delta.
\end{equation}
Let $\mathbf{R}_z = \mathbf{R}_z(V)$ be an arbitrary
generalized resolvent of $V$, and $C(\lambda;z_0)\in \mathcal{S}(N_{z_0};N_{ \frac{1}{ \overline{z_0} } })$
corresponds to $\mathbf{R}_z(V)$ by Inin's formula~(\ref{f1_7}).
$\mathbf{R}_z(V)$ has an analytic continuation to the set $\mathbb{D}\cup\mathbb{D}_e\cup \Delta$,
if and only if the following conditions are satisfied:

\begin{itemize}
\item[1)] $C(\lambda;z_0)$ has an extension to the set $\mathbb{D}\cup\Delta$ which is continuous in the
uniform operator topology;

\item[2)] The extended $C(\lambda;z_0)$ maps isometrically $N_{z_0}$ on the whole
$N_{ \frac{1}{ \overline{z_0} } }$, for all $\lambda\in\Delta$;

\item[3)] The operator $C(\lambda;z_0)- W_{\lambda;z_0}$ is invertible for all $\lambda\in\Delta$, and
\begin{equation}
\label{f5_29}
(C(\lambda;z_0)- W_{\lambda;z_0}) N_{z_0} = N_{ \frac{1}{ \overline{z_0} } },\qquad
\forall\lambda\in\Delta.
\end{equation}
Here $W_{\lambda;z_0}$ is defined by the following equality:
\begin{equation}
\label{f5_30}
W_{\lambda;z_0} P^H_{N_{z_0}} f = \frac{1-\overline{z_0}\lambda}{ \lambda - z_0}
P^H_{ N_{ \frac{1}{ \overline{z_0} } } } f,\qquad f\in N_\lambda,\ \lambda\in \mathbb{T}.
\end{equation}
\end{itemize}
\end{thm}
{\bf Proof. } We first notice that in the case~$z_0=0$ this theorem coincides with Theorem~\ref{t5_3}.
Thus, we may assume that $z_0\in \mathbb{D}\backslash\{ 0 \}$.

Let $\mathbf{R}_z(V)$ have an analytic continuation to the set $\mathbb{D}\cup\mathbb{D}_e\cup \Delta$.
Recall that the generalized resolvent $\mathbf{R}_z(V)$ is related to the generalized resolvent
$\mathbf{R}_z(V_{z_0})$ of $V_{z_0}$ by relation~(\ref{f2_2}), and that correspondence is bijective.
Then $\mathbf{R}_z(V_{z_0})$ admits an analytic continuation to $\mathbb{T}_e\cup \Delta_1$, where
$$ \Delta_1 = \left\{ \widetilde t:\ \widetilde t = \frac{\lambda - z_0}{1-\overline{z_0}\lambda},\
\lambda\in\Delta \right\}. $$
By Proposition~\ref{p2_1}, points $\widetilde t^{-1}$, $\widetilde t\in\Delta_1$, are of the regular type
of the operator $V_{z_0}$.
Moreover, relation~(\ref{f5_26}) written for the operator $V_{z_0}$ with $\zeta\in\Delta_1$, coincides with
relation~(\ref{f5_28}).
Then we may apply Theorem~\ref{t5_3} for the operator $V_{z_0}$ and the open arc $\Delta_1$.
If we then rewrite conditions 1)-3) of that theorem in terms of $C(\lambda;z_0)$, using the bijective
correspondence between $C(\lambda;z_0)$ for $V$, and $C(\lambda;0)$ for $V_{z_0}$,
we easily get conditions 1)-3) of the present theorem.

On the other hand, let conditions 1)-3) of the present theorem be  satisfied. Then conditions
of Theorem~\ref{t5_3} for the operator $V_{z_0}$ are satisfied. Therefore
$\mathbf{R}_z(V_{z_0})$ admits an analytic continuation to $\mathbb{T}_e\cup \Delta_1$.
Consequently, $\mathbf{R}_z(V)$ has an analytic continuation to the set $\mathbb{D}\cup\mathbb{D}_e\cup \Delta$.
$\Box$

\begin{center}
{\large\bf On the generalized resolvents for isometric operators with gaps.}
\end{center}
\begin{center}
{\bf S.M. Zagorodnyuk}
\end{center}

In this paper we obtain some slight correction and generalization of the results of Ryabtseva on
the generalized resolvents for isometric operators with a gap in their spectrum.
Also, analogs of some McKelvey's results and a short proof of Inin's formula for
the generalized resolvents of an isometric operator are obtained.

}
\end{document}